\documentclass[11pt,oneside,fleqn]{article}

\usepackage[ansinew]{inputenc}
\usepackage{verbatim}
\usepackage{amsmath,amssymb,amsthm}
\usepackage{graphicx}
\usepackage{multirow}
\usepackage{color}
\usepackage{cite}
\usepackage{subcaption}
\usepackage{hyperref}

\allowdisplaybreaks

\hypersetup{colorlinks, linkcolor=blue, citecolor=blue, urlcolor=blue}
\makeatletter
\long\def\@makecaption#1#2{%
  \vskip\abovecaptionskip\footnotesize
  \sbox\@tempboxa{#1. #2}%
  \ifdim \wd\@tempboxa >\hsize
    #1. #2\par
  \else
    \global \@minipagefalse
    \hb@xt@\hsize{\hfil\box\@tempboxa\hfil}%
  \fi
  \vskip\belowcaptionskip}
\makeatother

\newcommand{\todo}[1][\null]{\ensuremath{\clubsuit}}

\newcommand{\noprint}[1]{}

{\theoremstyle{definition}

\newtheorem*{remark*}{Remark}
}

\newcommand{\checked}[1][\null]{\ensuremath{\boldsymbol{\surd}}}

\begin{document}

\par\noindent {\LARGE\bf
Parallel stochastic methods for\\ PDE based grid generation
\par}

{\vspace{4mm}\par\noindent {\large Alexander Bihlo and Ronald D.\ Haynes
} \par\vspace{2mm}\par}

{\vspace{2mm}\par\noindent {\it
Department of Mathematics and Statistics, Memorial University of Newfoundland, \\
St.\ John's (NL), A1C 5S7, Canada\\
}}
{\noindent \vspace{2mm}{\it
\textup{E-mail}: abihlo@mun.ca, rhaynes@mun.ca
}\par}

\vspace{4mm}\par\noindent\hspace*{8mm}\parbox{140mm}{\small
The efficient generation of meshes is an important step in the numerical solution of various problems in physics and engineering.  We are interested in situations where global mesh
quality and tight coupling to the physical solution is important.  We  consider elliptic PDE
based mesh generation and present a method for the construction of adaptive meshes in two spatial dimensions using domain decomposition that is suitable for an implementation on parallel computing architectures. The method uses the stochastic representation of the exact solution of a linear mesh generator of Winslow type to find the points of the adaptive mesh along the subdomain interfaces. The meshes over the single subdomains can then be obtained completely
independently of each other using the probabilistically computed solutions along the interfaces as
boundary conditions for the linear mesh generator.
Further to the previously acknowledged performance characteristics, we demonstrate how
the stochastic domain decomposition approach is particularly suited to the problem of grid generation --- generating quality meshes efficiently.   In addition we show
further improvements are possible using interpolation of the subdomain interfaces and smoothing of mesh candidates.    An optimal placement strategy is introduced to automatically choose the number and placement of points along the interface using the mesh density function.
Various examples of meshes constructed using this stochastic--deterministic domain decomposition technique are shown and compared to the respective single domain solutions using a representative mesh quality
measure.
A brief performance study is included to show the viability of the stochastic domain decomposition approach and to illustrate the effect of algorithmic choices on the solver's efficiency.
}\par\vspace{2mm}

\section{Introduction}

The generation of flexible and adaptive meshes is an important problem in
several fields, such as computer graphics, the simulations of deformable
objects, engineering and scientific computing. See
e.g.~\cite{alli05a,frey10a,huan10a,shew02a,tera05a,thompson85,walk13a,yerr84a}
and references therein for an overview and applications to different
areas.

Our interest here is the use of grid adaptation as an
integral component of the numerical solution of many partial
differential equations (PDEs).  Here we are interested in the
calculation of an adaptive grid automatically tuned to the underlying
solution behaviour. The grid
is found by solving a mesh PDE which is (often) coupled
to the physical PDE of interest.
 Recently this general approach
has shown great promise, solving problems in
 meteorology \cite{budd2012},
relativistic magnetohydrodynamics \cite{he2012},
groundwater flow
\cite{Huang2005},
semiconductor devices \cite{yuan2012},
and viscoelastic flows \cite{zhang2011}, to name just a few.
A thorough, recent overview of PDE based mesh
generation may be found in \cite{huan10a}.  Such grid calculations
potentially add a burdensome overhead to the solution of the physical
model.    Here we present an efficient, parallel strategy for the solution of the mesh PDE based on a stochastic domain decomposition method recently proposed by Acebr{\'o}n et al.\ \cite{aceb05a}.
In many situations it is mesh quality, not an extremely accurate solution of the mesh PDE, which is
important.  It is for this reason the stochastic domain decomposition approach is a viable alternative
for mesh generation.

PDE based mesh generation may be classified by the type of PDE, elliptic or parabolic, linear or nonlinear, to be solved for the mesh.
The {\em equipotential} method of
mesh generation in 2D was first presented, to the best of our knowledge, by Crowley \cite{crowley62}.  The mesh
lines in the physical co-ordinates $x$ and $y$
are the level curves of the potentials $\xi$ and
$\eta$ satisfying  Laplace's equations
\begin{equation}\label{eq:crowley}\nabla^2 \xi  =0, \qquad \nabla^2 \eta = 0,\end{equation}
and  appropriate boundary conditions which ensure grid lines lie along the boundary of
the domain.
The mesh transformation,
giving the physical co-ordinates $x(\xi,\eta)$ and $y(\xi,\eta)$ in the physical domain
$\Omega_p$,
can be found either by (inverse) interpolation of the solution of (\ref{eq:crowley}) onto
 a (say) uniform $(\xi,\eta)$ grid or as Winslow \cite{wins66a} showed, by directly
solving the inverse equations of
(\ref{eq:crowley}), namely
\begin{subequations}\label{eq:winslow}
\begin{equation}\label{eq:winslowa}
\alpha x_{\xi\xi}-2\beta x_{\xi\eta}+\gamma x_{\eta\eta}=0,\quad
\alpha y_{\xi\xi}-2\beta y_{\xi\eta}+\gamma y_{\eta\eta}=0,
\end{equation}
where
\begin{equation}\label{eq:coeffs}
\alpha = x_\eta^2+y_\eta^2,\quad \beta=x_\xi x_\eta+y_\xi y_\eta
\quad\text{and}\quad \gamma = x_\xi^2+y_\xi^2.
\end{equation}
\end{subequations}
System (\ref{eq:winslow}) can be solved directly for the mesh
transformations $x(\xi,\eta)$ and $y(\xi,\eta)$ using a uniform
grid for the variables $\xi$ and $\eta$ belonging to the artificial
computational domain $\Omega_c$.

In this paper, we will exclusively work with mesh generators defined on the physical space, i.e.\ yielding $\xi=\xi(x,y)$ and $\eta=\eta(x,y)$. As mentioned, a convenient way to invert the mesh transformation from the physical space $\Omega_{\rm p}$ to the computational space $\Omega_{\rm c}$ is through \emph{interpolation}. That is, once the numerical solution yielding $\xi_{ij}=\xi(x_i,y_j)$, $\eta_{ij}=\eta(x_i,y_j)$ is obtained one obtains the values for $x_{ij}=x(\xi_i,\eta_j)$, $y_{ij}=y(\xi_i,\eta_j)$ from two-dimensional interpolation. An alternative would be to numerically solve the hodograph transformations
\[
 \left(\begin{array}{cc}x_\xi & x_\eta\\ y_\xi & y_\eta\end{array}\right)=\frac{1}{J}\left(\begin{array}{cc} \eta_y & -\xi_y\\ -\eta_x & \xi_x\end{array} \right),
\]
where $J=\xi_x\eta_y-\xi_y\eta_x$. As this leads to another system of PDEs,  we choose the interpolation method to obtain the physical mesh lines. In practice, this inversion to the physical co--ordinates is not necessary. Instead one could transform the physical PDE of interest into the computational co--ordinate system.

Dvinksy \cite{dvin91a} discusses the existence and uniqueness
of such mesh transformations in the context of one--to--one harmonic maps.
He demonstrates that the solutions to (\ref{eq:crowley}), and hence (\ref{eq:winslow}),
will be well--defined if the $(\xi,\eta)$ domain is convex.  Since we construct
this domain ourselves this condition can always be satisfied.

Godunov and Prokopov \cite{godunov72},  Thompson et al.\
\cite{thompson74,thompson85} and Anderson \cite{anderson87}, for example, add terms to
(\ref{eq:crowley})
 and (\ref{eq:winslow}) to better control the mesh distribution and
quality.
Winslow \cite{wins81a} generalizes (\ref{eq:crowley})
by adding a diffusion coefficient $w(x,y) > 0 $ depending on the
gradient or other aspects of the solution.  This gives
the linear elliptic mesh generator
\begin{equation}\label{eq:linearmesh}
\nabla\cdot(w\nabla \xi) = 0\quad\text{and}\quad \nabla\cdot(w\nabla \eta) = 0,
\end{equation}
of interest in this paper.  The diffusion
coefficient $w$ characterizes regions where additional
mesh resolution is needed.
An example for a mesh in both computational and physical co--ordinates is depicted in Figure~\ref{fig:SingleDomainSolutionRunningExample} on page \pageref{fig:SingleDomainSolutionRunningExample}.  This figure shows that the grids in the computational and physical co--ordinates behave inversely; the areas of grid concentration of the mesh in the computational co--ordinates coincide with the areas of de-concentration in the physical space and vice versa.

There has been recent work to parallelize nonlinear PDE based mesh generation using  on Schwarz based domain decomposition
approaches.  In \cite{Haynes:2012},
Haynes and Gander propose and analyze classical, optimal and optimized
Schwarz methods in one spatial dimension.   A numerical study of
classical and optimized Schwarz domain decomposition for 2D nonlinear mesh generation
has been presented in \cite{hayn12a}.
Recently, in \cite{Chen:2012uo}, a  {\em monolithic} domain decomposition method, simultaneously solving a mesh generator similar to (\ref{eq:WinslowTypeGridGenerator}) coupled to the physical PDE, was presented for a shape
optimization problem.  The authors used an overlapping domain decomposition approach
to solve the coupled problem.  Here we focus on the linear mesh generation problem only, detailing the effect of our stochastic domain decomposition method on the generated meshes.

This paper is organized as follows. In Section~\ref{sec:TheoryDD} we review the necessary background material on the stochastic interpretation of solutions of linear elliptic boundary value problems and their relation to numerical grid generators. We also briefly discuss the techniques required to achieve
such stochastic solutions numerically. We explain how to couple
the stochastic solution of linear elliptic mesh generators with a domain decomposition approach
to obtain a scalable version of the algorithm.
The parallel performance of the approach is also reviewed.
Section \ref{sec:DDStochastic}
illustrates our parallel grid generation strategy.
The effect of smoothing the probabilistically computed interface solutions and the subdomain
solutions is  demonstrated.
Section~\ref{sec:ResultsDD} is devoted to further examples of grids computed using stochastic domain
decomposition.
In Section~\ref{sec:PerformanceStudy} we illustrate the performance of our stochastic domain decomposition
method for grid generation as compared to the single domain, global
solution strategy, paying particular attention to the effect of our algorithmic choices.
Section~\ref{sec:ConclusionsDD} contains the conclusions of the paper as well as thoughts for further research directions.

\section{Winslow mesh generation using a stochastic domain decomposition method}\label{sec:TheoryDD}

In this section, following \cite{aceb05a},  we describe how to generate adaptive meshes
by solving Winslow's mesh generator~(\ref{eq:linearmesh}) using a stochastic representation
of the solution along the artificial interfaces and a non--overlapping, non--iterative domain decomposition.

\subsection{Background}

For concreteness, we will consider a two-dimensional linear grid generator of Winslow type
\begin{align}\label{eq:WinslowTypeGridGenerator}
\begin{split}
 &-\nabla\cdot\left(w\nabla\xi\right)=0,\quad -\nabla\cdot\left(w\nabla\eta\right)=0,\\
 &\xi|_{\partial \Omega_{\rm p}}=f(x,y),\quad \eta|_{\partial \Omega_{\rm p}}=g(x,y),
\end{split}
\end{align}
where $w=w(x,y)>0$ is a strictly positive weight function~\cite{huan10a,wins66a}. This grid generator
finds  the stationary solution of spatially dependent diffusion processes, yielding the computational coordinates $\xi=\xi(x,y)$ and $\eta=\eta(x,y)$ in terms of the physical coordinates $x$ and $y$.
For the sake of simplicity, we restrict ourselves to the case of
rectangular physical and computational domains.   As mentioned, Dvinksy \cite{dvin91a}, then guarantees
a well--posed mesh generation problem.


Equations~\eqref{eq:WinslowTypeGridGenerator} form a system of two decoupled linear elliptic PDEs. It is well known that for such boundary value problems the solution can be written using methods of stochastic calculus~\cite{kara91a,okse10a}. In the present case, this solution is conveniently derived from the expanded form of system~\eqref{eq:WinslowTypeGridGenerator}, which reads
\begin{equation}\label{eq:WinslowTypeGridGeneratorExpanded}
 \frac1w\nabla w\cdot \nabla\xi+\nabla^2\xi=0,\quad \frac1w\nabla w\cdot \nabla\eta+\nabla^2\eta=0.
\end{equation}
The solution of system~\eqref{eq:WinslowTypeGridGenerator} has the stochastic representation
\begin{subequations}\label{eq:WinslowTypeGridGeneratorStochasticSolution}
\begin{equation}\label{eq:WinslowTypeGridGeneratorStochasticSolutionA}
 \xi(x,y)=\mathrm{E}[f(\mathbf{X}(\tau))],\quad \eta(x,y)=\mathrm{E}[g(\mathbf{X}(\tau))],
\end{equation}
where $\mathbf{X}(t)=(x(t),y(t))^{\rm T}$ satisfies, in the \^{I}to sense, the stochastic differential equation
\begin{equation}\label{eq:WinslowTypeGridGeneratorStochasticSolutionB}
 \mathrm{d}\mathbf{X}(t)=\frac1w\nabla w\,\mathrm{d}t+\mathrm{d}\mathbf{W}(t).
\end{equation}
\end{subequations}
In the stochastic solution~\eqref{eq:WinslowTypeGridGeneratorStochasticSolution}, $\mathrm{E}[\cdot]$ denotes the expected value, $\tau$ is the time when the stochastic path starting at the point $(x,y)$ first hits the boundary
of the physical domain $\Omega_{\rm p}$ and $\mathbf{W}$ is the standard two-dimensional Brownian motion~\cite{aceb05a,kara91a,okse10a}.

As noted in~\cite{aceb05a}, seeking the numerical solution of system~\eqref{eq:WinslowTypeGridGenerator} using the probabilistic solution at \emph{all} the grid points of interest is generally too expensive, especially when compared to direct or iterative (deterministic) methods
to solve the linear system of equations resulting from the discretization of (\ref{eq:WinslowTypeGridGenerator}).
Efficiency is indeed a crucial factor for grid adaptation, since the grid
is computed in addition to the solution of the physical PDE.  This is particularly important
for time-dependent problems where~\eqref{eq:WinslowTypeGridGenerator} has to be solved in combination with a system of PDEs at every time step.

\begin{figure}[!ht]
\centering
 \includegraphics[scale=1]{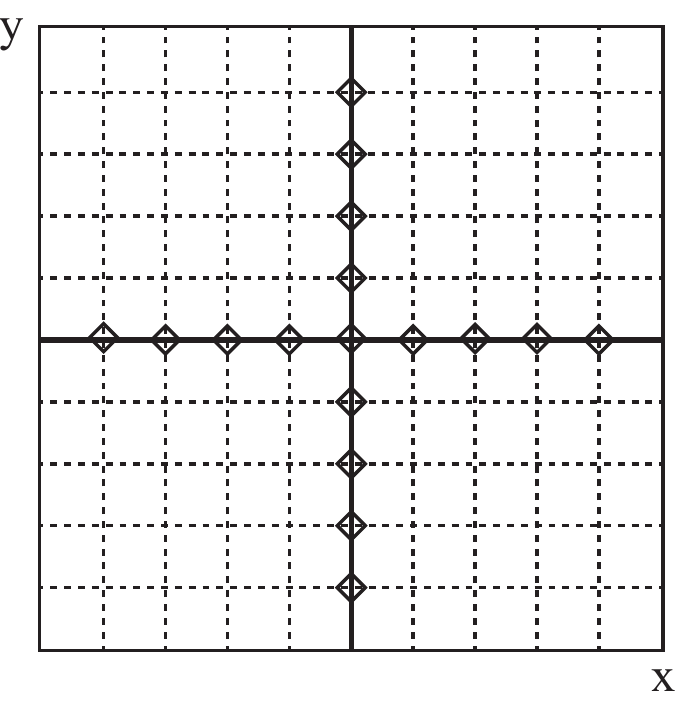}
\caption{Domain decomposition using the stochastic approach. The probabilistic solution~\eqref{eq:WinslowTypeGridGeneratorStochasticSolution} of the linear mesh generator~\eqref{eq:WinslowTypeGridGenerator} is used to compute the exact solution $\xi(x,y)$, $\eta(x,y)$ at the interface points (open circles). The solution on each subdomain (the dashed lines) is computed in parallel using a suitable single domain solver with the interface solutions serving as Dirichlet boundary values.}
\label{fig:CartoonOnProbabilisticDD}
\end{figure}

To address the expense of the stochastic approach, the key idea put forward in~\cite{aceb05a} is to use the probabilistic solution~\eqref{eq:WinslowTypeGridGeneratorStochasticSolution} in the context of domain decomposition.  The probabilistic solution is used
only to obtain the boundary conditions at the subdomain interfaces, see
Figure~\ref{fig:CartoonOnProbabilisticDD}.

This approach to domain decomposition has the advantage of being fully parallelizable, as each subdomain can be assigned to a single core and thus the grids at the various subdomains can be computed completely independently of each other. Once the interface values are
obtained with sufficient accuracy the solutions on the subdomains are only computed once --- no
iteration is needed in this domain decomposition approach.  In addition, the two-dimensional random walks required to solve the stochastic differential equations~\eqref{eq:WinslowTypeGridGeneratorStochasticSolutionB} at each interface point can be done
independently and efficiently using, for example, an implementation on a graphics processing unit (GPU).  We stress here that several Monte Carlo techniques have already been successfully realized on GPUs, e.g.\ for problems in molecular dynamics, statistical physics and financial mathematics (see~\cite{prei09a,vanM08a} for examples), leading to a massive speed-up (often several orders of magnitude) compared to the conventional CPU solution. The perspective of compatibility with the principles of GPU programming is therefore a highly desirable feature of a new numerical algorithm.

In theory, the extension of Winslow type mesh generator (\ref{eq:WinslowTypeGridGenerator}) and the
stochastic representation of the solution (\ref{eq:WinslowTypeGridGeneratorStochasticSolutionA}) and
(\ref{eq:WinslowTypeGridGeneratorStochasticSolutionB}) are readily extensible to three spatial dimensions.  Indeed, the efficacy of the Monte Carlo approach is known to increase in higher dimensions.
Realistically, there is always a challenge to move to a three dimensional implementation.
For instance, one must select points along, the now two dimensional artificial interface,
to evaluate the stochastic representation of the solution.   Two dimensional interpolation
is then required to provide the boundary conditions for the subdomain solves.


\subsection{Implementation details}

We now describe the solution procedure for solving the system of stochastic differential equations~\eqref{eq:WinslowTypeGridGeneratorStochasticSolutionB}. This system of equations is particularly simple in that it is decoupled, i.e.\ the equations for $x(t)$ and $y(t)$ can be solved independently of each other and hence this solution is fully parallelizable. The components of the two-dimensional Brownian motion $\mathrm{d} \mathbf{W}$ can be realized as $\mathrm{d}W_i=2\sqrt{\Delta t}\, N(0,1)$, for $i=1,2$, where $N(0,1)$ is a normally distributed random number with zero mean and variance one~\cite{kara91a}. A simple method for integrating the stochastic differential equations~\eqref{eq:WinslowTypeGridGeneratorStochasticSolutionB} is to use the Euler--Maruyama method~\cite{kloe92a} with a time step $\Delta t$.

The crucial part of solving~\eqref{eq:WinslowTypeGridGeneratorStochasticSolution} is the correct determination of the \textit{exit time} $\tau$, i.e.\ the time when the stochastic process starting at $(x,y)$ first hits one of the boundaries of $\Omega_p$. One of the simplest possibilities to estimate the exit time is to integrate~\eqref{eq:WinslowTypeGridGeneratorStochasticSolutionB} using the Euler--Maruyama method and to determine whether the process leaves the domain at the end of the time step. The problem with this approach is that one cannot determine whether the true stochastic process already exited the domain during the time step. It is possible that the true process started within the domain, left it, but returned to the domain before the end of the time step.
Indeed, we found that the plain Euler--Maruyama method leads to larger errors in the fully probabilistically generated grid lines close to the boundaries. The situation is much approved if a linear Brownian bridge is used to check if the stochastic process left the domain within a time-step. Using this Brownian bridge as an interpolating process, a simple expression for the conditional exit probability of the process can be found, see~\cite{gobe00a} for further details. Using this technique greatly improves the performance of linear time-stepping.

An alternative to linear time-stepping is to implement the approach proposed in~\cite{jans03a}, based on \textit{exponential time-stepping}.
Unlike the Euler--Maruyama scheme, which uses a constant time step $\Delta t$, this approaches
chooses $\Delta t$ as an
exponentially distributed random variable. That is, each time step is determined as an independent realization of a distribution with the probability density function $\lambda\exp(-\lambda t)$, where $\lambda>0$ is a constant parameter. The expected value of an exponentially distributed random variable $X$
is then $\mathrm{E}[X]=1/\lambda$, i.e.\ as $\lambda\to\infty$, $\Delta t\to0$.

The central idea behind exponential time-stepping is that it is possible to carry out the boundary test for the exit time by both explicitly checking if the process has exited the domain and also using the conditional probability that the process reached a boundary between $\mathbf{X}(t)$ and $\mathbf{X}(t+\delta t)$~\cite{jans03a}. Exponential time-stepping works because this conditional probability is a function of the boundary itself. It was shown in~\cite{jans03a} that for several classes of stochastic differential equations this conditional probability can be either computed analytically or approximated numerically with high efficiency.  System~\eqref{eq:WinslowTypeGridGeneratorStochasticSolutionB} constitutes a process of the \emph{additive noise} class --- a deterministic system (the first term in~\eqref{eq:WinslowTypeGridGeneratorStochasticSolutionB}) that is superimposed by constant noise (the second term in~\eqref{eq:WinslowTypeGridGeneratorStochasticSolutionB}). For such classes of problems, the precise implementation of the exponential time-stepping method is presented in Section 4.2 of~\cite{jans03a}.

Below we present examples for both time-stepping techniques. We found that as long as a boundary test is applied in linear time-stepping, both methods for solving the SDE lead to grids that are sufficiently accurate near the boundaries.

The system of stochastic differential equations~\eqref{eq:WinslowTypeGridGeneratorStochasticSolutionB} is solved many times, as one needs enough samples to approximate the expected value in~\eqref{eq:WinslowTypeGridGeneratorStochasticSolutionA} through the arithmetic mean.
Indeed the main problem with Monte Carlo techniques is the rather slow convergence rate.
The error in estimating the expected value with the sample mean, using pseudo--random
numbers, is known to be of the order of $N^{-1/2}$, where $N$ is the number of samples,
see e.g.~\cite{pres07a}.

The main problem with
pseudo-random numbers is that a significant fraction of them tend to cluster in localized regions of the sampling space,
while other parts of the domain are under-sampled. This variation in the point density is measured through the \emph{discrepancy}.   Sequences of numbers that fill the sampling space while
reducing the discrepancy are called quasi-random numbers. They can lead to a significant speed-up of the Monte Carlo method, yielding convergence rates of order $N^{-1}$~\cite{pres07a}. The downside of quasi-random numbers is that they cannot be generated in a fully parallel way. Due to their inherent correlation it is necessary to first generate a sequence and then re-order the single elements of this sequence to break the correlation amongst the elements. This procedure is not a parallel operation.

Despite the superior convergence properties of Monte Carlo techniques employing quasi-random numbers,
 we choose to use pseudo--random numbers in what follows.   Given the
 relative low accuracy needed when solving the mesh PDE, we find that
 the use of a moderate number of Monte Carlo simulations ($N\approx 10^{4}$) with
 pseudo--random numbers gives sufficiently accurate  meshes. Moreover, as shown in the
 next section, applying a smoothing procedure to the computed mesh or simply using Monte Carlo solutions only along the subdomain interfaces has the potential to further reduce the number of Monte Carlo simulations needed to obtain a good quality mesh.
 At the same time,  avoiding the quasi-random numbers prevents the introduction of a bottleneck in the otherwise fully parallel algorithm of stochastic domain decomposition for grid
 generation.

\subsection{Parallel performance}

An analysis of the parallel performance of this stochastic domain decomposition (SDD) algorithm was provided in \cite{aceb05a}.  We summarize their findings here.  In two spatial dimensions, a straightforward calculation shows that the speedup $S_p$ obtained by using $p$ processors (equal to the number of subdomains $(n_x+1)(n_y+1)$) is
\begin{equation}\label{eq:speedup2}
S^{\rm SDD}_p := \frac{p}{1+k(n_x+n_y)\frac{T_{\rm MC}}{T_1}}.
\end{equation}
In (\ref{eq:speedup2}) it is assumed that we compute the solution stochastically at $k$ points along each of the $(n_x+n_y)$  interior
interfaces.  The quantity $T_{\rm MC}$ is the time to compute the solution at a single interface point using
Monte Carlo simulations, while $T_1$ is the time to compute the solution sequentially on the entire domain.
The total time spent on the stochastic portion of the algorithm, $T_{\rm stoc}$, is simply
$T_{\rm stoc}= k(n_x+n_y)T_{\rm MC}$.
This formula is arrived at by modelling the compute time on $p$ processors by
$$T_p = \frac{k(n_x+n_y)T_{\rm MC}}{p}+\frac{T_1}{p},$$
which assumes the subdomain solves take no more than $1/p$ of the global solve time $T_1$ and
we split the time for the $k(n_x+n_y)$ Monte Carlo simulations amongst the $p$ processors.
The speed--up is then found simply as $S_p = T_1/T_p$.
If we assume $n_x=n_y=:n$ then as $n\rightarrow \infty$ we obtain
$$S^{\rm SDD}_p \sim \frac{T_1}{2kT_{\rm MC}}\sqrt{p}.$$
In $d$--dimensions this generalizes to
$$S^{\rm SDD}_p \sim \frac{T_1}{dkT_{\rm MC}}p^{1-1/d}.$$

For parallel finite differences (PFD) using a deterministic domain decomposition approach we have
\begin{equation}\label{eq:speeduppfd}
S^{PFD}_p:=\frac{p}{1+\frac{T_{\rm com}}{T_1}p}\sim \frac{T_1}{T_{\rm com}} \quad \text{as}\,\,p\rightarrow \infty,
\end{equation}
where $T_{\rm com}$ is the time spent on communications between the processors.  Compared to
(\ref{eq:speedup2}) we see how the PFD approach becomes saturated; scaling is affected as $p\rightarrow\infty$.

Acebr\'{o}n et al.\ \cite{aceb05a} proceed to model $T_{\rm MC}, T_1$ and $T_{\rm com}$ and conclude that
the speedup of the parallel stochastic algorithm is better than that of the parallel finite difference (PFD) method if
the number of available processors is large enough.  Moreover, if the spatial dimension is large enough then
$S^{\rm SDD}_p > S^{\rm PFD}_p$ even for small numbers of processors.

In practice, the numerics shown in \cite{aceb05a} show that the SDD approach outperforms PFD in terms of total computation time with the difference growing as the number of processors increases.  As expected, for
small numbers of processors the majority of the computation time need for SDD is taken with the local (deterministic) subdomain
solves.   As the number of processors increases proportionately more time is needed for the Monte Carlo
simulations.  In \cite{acebron:2009fi,acebron:2010ce} and \cite{Acebron:2011kd}, the authors also compare the SDD approach to the
widely used parallel ScaLAPACK library which is used to find the whole domain solution in parallel.
In these papers, the authors extend the technique to handle nonlinear problems for which a
probabilistic representation exists.
The examples shown clearly demonstrate the better
performance of SDD compared to widely available parallel alternatives.

In what follows, we focus on the merits of the SDD approach for PDE based
mesh generation.  In the original papers~\cite{aceb05a,acebron:2010ce} it was extensively shown that SDD is suitable for numerically solving certain classes of physical PDEs for which a
stochastic representation of the solution is known.   In PDE based mesh generation the main task is not to solve the mesh PDEs with high accuracy but to obtain high quality meshes. This leads to both
new challenges and opportunities for additional speedup.
In particular, we propose an automatic strategy to select the number and location of the points along the interface at which the stochastic representation will be evaluated.
As adaptive meshes generally evolve with the solution of the physical PDEs,  fixed choices of the number and position of these points along the interface are generally not advisable.
This dynamic choice has the potential to increase the efficiency of the SDD approach.
Furthermore we show that it is possible to reduce the number of Monte Carlo simulations at each chosen interface point and still obtain quality meshes by appropriately smoothing the mesh solutions.
Although Acebr\'{o}n et al.\ have clearly demonstrated the parallel efficiency of the SDD approach, we include a small timing study to corroborate these results in the face of the specific algorithmic choices
we make for the PDE based mesh generation problem.

\section{A numerical case study}\label{sec:DDStochastic}

In this section we illustrate the approach described in the previous section, taking care to
detail the effect of each algorithmic choice on the resulting mesh.
We illustrate the successive stages of our implementation by
 solving the Winslow mesh generator (\ref{eq:WinslowTypeGridGenerator}) with $w=1/\rho$ where $\rho$ is
the monitor function
\begin{equation}\label{eq:MonitorFunctionRunningExample}
 \rho=1+R\exp(-50(x-3/4)^2-50(y-1/2)^2-1/2).
\end{equation}
In all examples we choose $R=15$.
In the physical co-ordinates, this monitor function will concentrate the  mesh
near $x=3/4$ and $y=1/2$.
It is instructive to note that due to the relation $w=1/\rho$, the region where~$\rho$ attains its maximal values is where the mesh in the \emph{physical space} is the sparsest.
We begin by showing the mesh obtained by solving (\ref{eq:WinslowTypeGridGenerator}) on
single domain, displaying the mesh in both computational and physical co-ordinates.   The entire mesh
is then recomputed using the stochastic method.  The effect of the number of Monte Carlo
simulations is shown.   We then compute the mesh stochastically only along the artificial
interfaces; the rest of the mesh is computed using a domain decomposition approach with deterministic subdomain solves.
We finally show how smoothing can be incorporated to obtain quality
meshes while keeping the number of Monte Carlo simulations small.

\subsection{Mesh quality measures}

The effects of the different steps of the proposed algorithm on the meshes
can be quantified by introducing a  \emph{mesh quality measure}.  Several mesh quality measures have been proposed in the literature to quantify properties of adaptive meshes. These measures usually assess mesh regularity, the degree of adaptivity to the numerical solution, the regularity of mesh elements in an appropriate metric (defined through the monitor function) or they quantify equidistribution~\cite{huan10a}.

In the present case, we restrict ourselves to a comparison of the mesh quality of the domain decomposition solution with the mesh quality of the single domain mesh. We do not aim to
evaluate the absolute mesh quality of the single domain solution of the linear mesh generator~\eqref{eq:WinslowTypeGridGenerator} here,  but rather
we are interested in estimating how well the domain decomposition solution approximates the quality of the single domain grid. This is an important task because a mesh quality measure can yield a \emph{stopping criterion} for the probabilistic algorithm as we are not concerned with finding a perfect numerical solution of the mesh generator~\eqref{eq:WinslowTypeGridGenerator}.


For the sake of simplicity,  we will exclusively work with the \emph{geometric mesh quality measure} defined by
\begin{equation}\label{eq:Qmeasure}
 Q(K)=\frac12\frac{\mathrm{tr}(J^{\rm T}J)}{\sqrt{\det(J^{\rm T}J)}},
\end{equation}
where $J$ is the Jacobian of the transformation $x=x(\xi,\eta)$, $y=y(\xi,\eta)$ and $K$ is a mesh element in $\Omega_{\rm c}$, see~\cite{huan10a} for more details. This measure has the property $Q(K)\ge1$ with $Q(K)=1$ for an equilateral mesh element only.  Probabilistically computed meshes that have not yet converged usually show several kinks in mesh lines and thus feature grid cells that are far away from equilateral
as compared to the meshes obtained from a deterministic algorithm.  Hence $Q(K)$
is an appropriate measure for our purposes.

\subsection{Single domain reference mesh}

We first solve system~\eqref{eq:WinslowTypeGridGenerator} with monitor function~\eqref{eq:MonitorFunctionRunningExample} on the entire domain $\Omega_{\rm p}=[0,1]\times[0,1]$.
This mesh will serve as the reference solution for the stochastically generated
meshes which follow.


System~\eqref{eq:WinslowTypeGridGenerator} is discretized with centered finite differences on
a uniformly spaced grid in the physical $(x,y)$ co--ordinates and the resulting system is solved using a Jacobi iteration. We use the Dirichlet boundary conditions $\xi(0,y)=0$, $\xi(1,y)=1$, $\eta(x,0)=0$ and $\eta(x,1)=1$.  This ensures grid lines on
the boundary of the unit square.  The values of $\xi(x,0)$, $\xi(x,1)$, $\eta(0,y)$ and $\eta(1,y)$ are obtained from solving the one-dimensional forms of the mesh generator~\eqref{eq:WinslowTypeGridGenerator}.

\begin{figure}[!ht]
        \centering
        \begin{subfigure}[b]{0.45\textwidth}
                \centering
                \includegraphics[width=\textwidth]{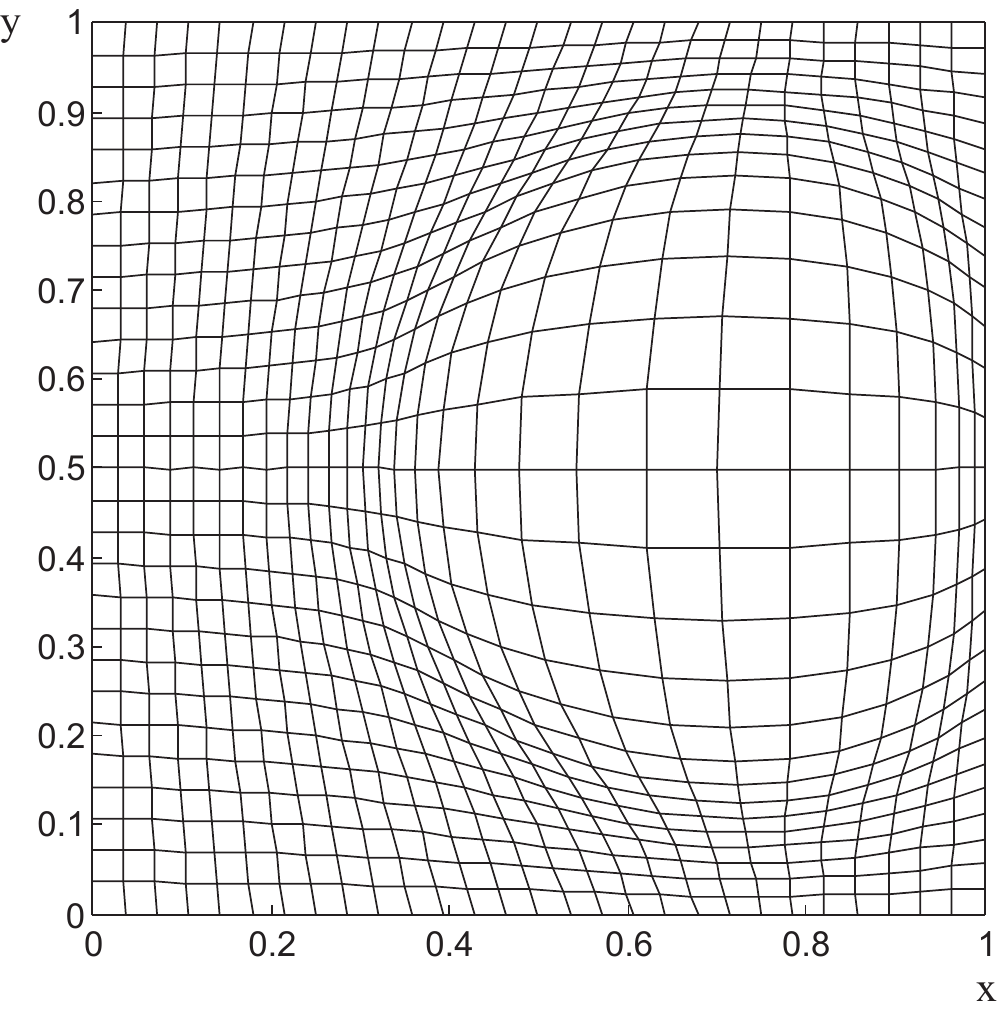}
        \end{subfigure}\qquad
        \begin{subfigure}[b]{0.45\textwidth}
                \centering
                \includegraphics[width=\textwidth]{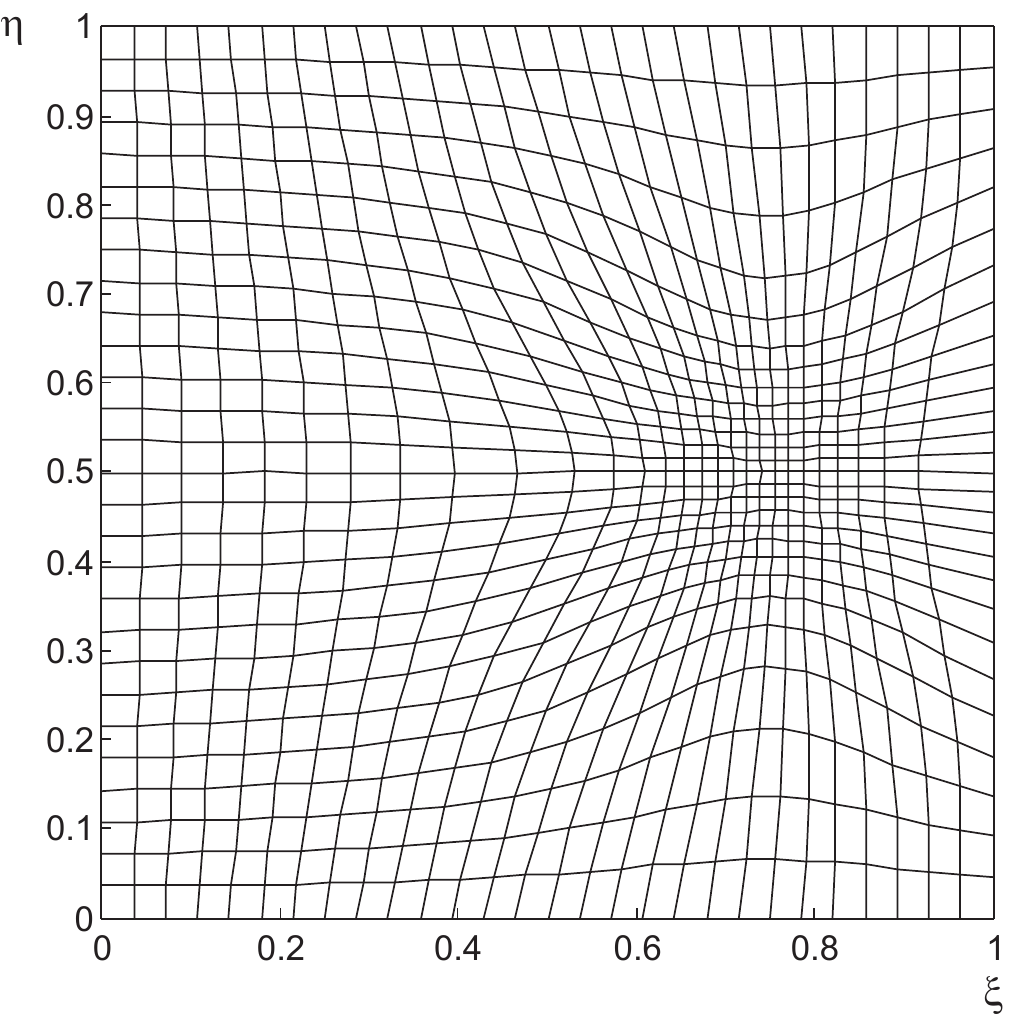}
        \end{subfigure}
        \caption{\textbf{Left:} Numerical solution of the Winslow-type grid generator~\eqref{eq:WinslowTypeGridGenerator} for the monitor function~\eqref{eq:MonitorFunctionRunningExample} over the physical space $\Omega_{\rm p}=[0,1]\times [0,1]$. \textbf{Right:} The physical mesh lines obtained from the grid
        on the left using linear interpolation.}
\label{fig:SingleDomainSolutionRunningExample}
\end{figure}

Figure~\ref{fig:SingleDomainSolutionRunningExample} (left) shows the numerical approximation to
$\xi(x,y)$ and $\eta(x,y)$ obtained by solving~\eqref{eq:WinslowTypeGridGenerator} using the monitor function~\eqref{eq:MonitorFunctionRunningExample} on $\Omega_{\rm p}=[0,1]\times [0,1]$ with $29\times 29$ uniformly chosen grid points. On the right of Figure~\ref{fig:SingleDomainSolutionRunningExample} we depict the physical mesh lines obtained by interpolation of the mesh on the left onto uniformly
spaced $(\xi,\eta)$ co--ordinates. Here, as expected, the regions of mesh concentration indeed coincide with the maximal values of~$\rho$, near $x=3/4$ and $y=1/2$.

\subsection{Probabilistically computed mesh}

\begin{figure}[!ht]
\centering
\begin{subfigure}[b]{0.45\textwidth}
  \centering
  \includegraphics[width=\linewidth]{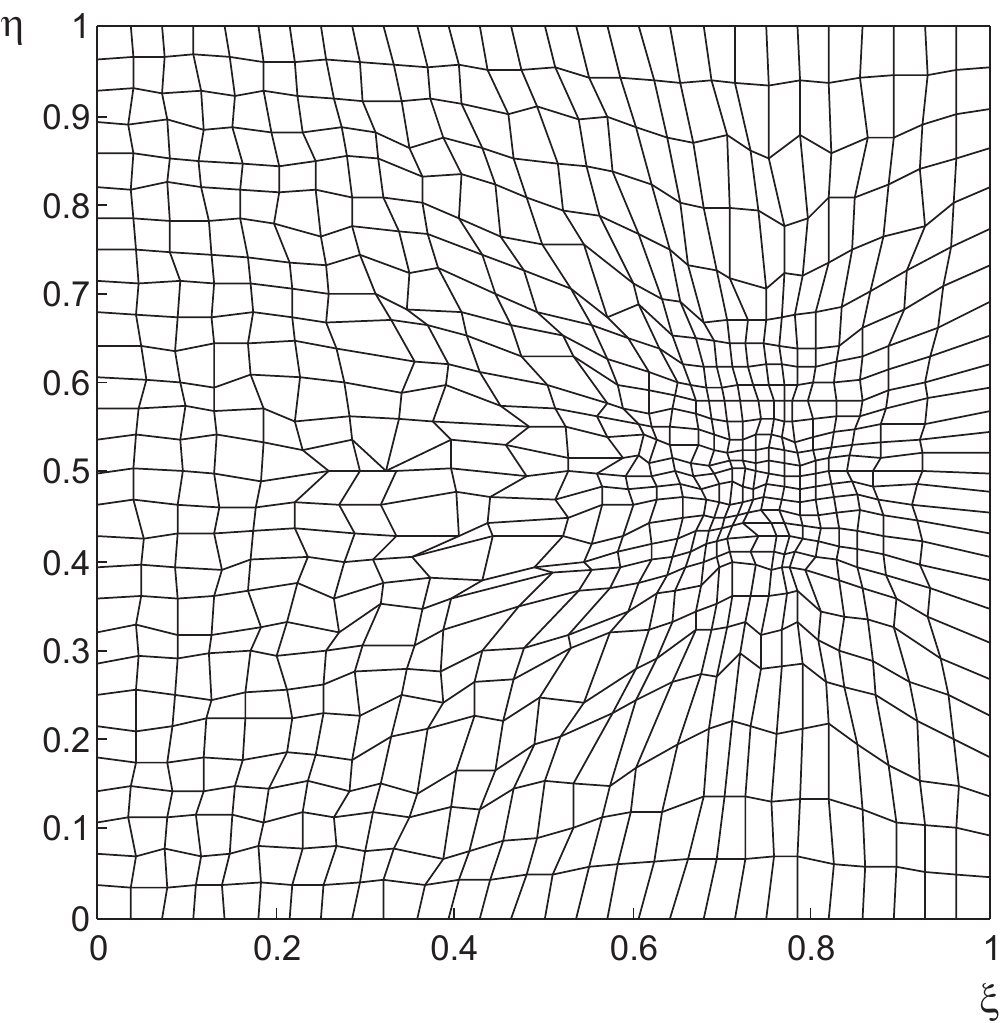}
  \label{fig:SingleDomainSolutionFullyProbabilisticRunningExample1000}
\end{subfigure}\qquad
\begin{subfigure}[b]{0.45\textwidth}
  \centering
  \includegraphics[width=\linewidth]{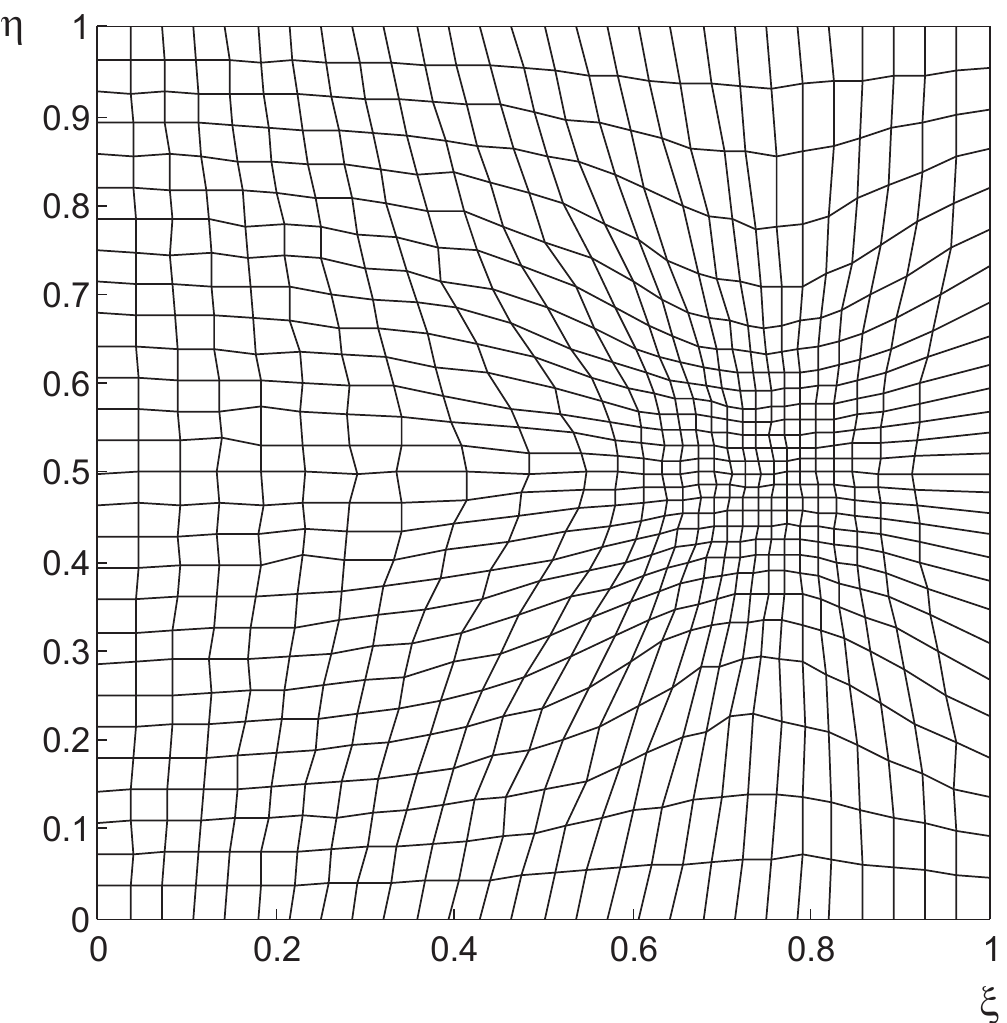}
  \label{fig:SingleDomainSolutionFullyProbabilisticRunningExample10000}
\end{subfigure}
\caption{\textbf{Left:} Grid obtained from the probabilistic solution~\eqref{eq:WinslowTypeGridGeneratorStochasticSolution} using $\lambda=1000$ and $N=1000$ Monte Carlo simulations at each point to estimate the expected value. \textbf{Right:} Same as left, but using $N=10000$ Monte Carlo simulations.}
\label{fig:SingleDomainSolutionFullyProbabilisticRunningExample}
\end{figure}

It is instructive to display the result that is achieved by solving the Winslow-type mesh generator~\eqref{eq:WinslowTypeGridGenerator} for \emph{all} the grid points using the stochastic solution~\eqref{eq:WinslowTypeGridGeneratorStochasticSolution}. The result of such fully probabilistically computed solutions is displayed in Figure~\ref{fig:SingleDomainSolutionFullyProbabilisticRunningExample}.

We choose the exponential time stepping parameter $\lambda=1000$ for both runs. The difference between Figure~\ref{fig:SingleDomainSolutionFullyProbabilisticRunningExample} (left) and Figure~\ref{fig:SingleDomainSolutionFullyProbabilisticRunningExample} (right) is the number of Monte Carlo simulations used.
On the left  only $N=1000$  Monte Carlo simulations were used to approximate the expected values in the solution~\eqref{eq:WinslowTypeGridGeneratorStochasticSolutionA}.  On the right we use $N=10000$ Monte Carlo runs.   Comparing to the reference solution on the right of Figure \ref{fig:SingleDomainSolutionRunningExample},
the difference in the results is dramatic.   This illustrates the well--known problem that a high number of simulations is needed to obtain the solution of a PDE with sufficient accuracy from a stochastic representation. As discussed in Section~\ref{sec:TheoryDD}, the theoretical explanation for this finding is the slow convergence rate of Monte Carlo methods, which is proportional to $N^{-1/2}$ if pseudo--random number are used during the Monte Carlo simulations. It is mainly this slow convergence which has prevented the wider use of Monte Carlo techniques in the numerical solution of PDEs.  Here we suggest that the relative lower accuracy requirement for the mesh makes the stochastic formulation, coupled with the domain decomposition
 approach to follow, an appealing strategy.

To quantify the above findings, we compute the geometric mesh quality measure (\ref{eq:Qmeasure}) for both the deterministic mesh and the two meshes
found by the stochastic algorithm. As $Q=Q(K)$ is a function of each mesh element $K$, we only list the maximum and the mean of $Q$ over the domain, denoted by $Q_{\rm max}$ and $Q_{\rm mean}$ for the deterministic solution and by $Q^{\rm S}_{\rm max}$ and $Q^{\rm S}_{\rm mean}$ for the stochastic solution.   Comparisons  of the stochastic solution and the deterministic solution are given by the ratios $R_{\rm max}=Q_{\rm max}/Q^{\rm S}_{\rm max}$ and $R_{\rm mean}=Q_{\rm mean}/Q^{\rm S}_{\rm mean}$.
We also present $l_\infty$, a measure of how well $\rho$ is captured by
the generated meshes.
More precisely, we compute the linear interpolant of $\rho$ onto the grids obtained using the deterministic and the stochastic solver yielding $\rho_{\rm interp}$ and $\rho_{\rm interp}^{\rm S}$, respectively. The $l_\infty$ error is then computed as the maximal difference between $\rho_{\rm interp}$ and $\rho_{\rm interp}^{\rm S}$, i.e.\ $l_\infty=\max(|\rho_{\rm interp}-\rho_{\rm interp}^{\rm S}|_i)$, where $i$ runs through all grid points. We have chosen this error measure as it allows us to
compare the maximal deviation of the interpolated $\rho$ on the stochastically
generated mesh to interpolated $\rho$ on the deterministically generated mesh,  and hence provides an additional measure of convergence for the Monte Carlo simulation.
Smaller values of $l_\infty$ suggest a better resolution of the mesh density function $\rho$.

\begin{table}[!ht]
  \centering
  \caption{Mesh quality for grid adapting to monitor function~\eqref{eq:MonitorFunctionRunningExample}. $\lambda=1000$.}
  \begin{tabular}{|p{4cm}|r|r|r|r|r|r|}
  \hline
   & $N$ & $l_{\infty}$ & $Q^{\rm S}_{\rm max}$ & $Q^{\rm S}_{\rm mean}$ & $R_{\rm max}$ & $R_{\rm mean}$ \\
   \hline
   \multirow{2}{*}{\parbox{4cm}{$Q_{\rm max}=1.8$\\ $Q_{\rm mean}=1.16$}}& $1000$ & $2.23$ & $3.85$ & $1.22$ & $0.47$ & $0.95$ \\
   &$10000$ & $1.33$ & $1.9$ & $1.18$ & $0.95$ & $0.98$\\
   \hline
   \end{tabular}
   \label{tab:MeshQualityMonteCarlo}
\end{table}

It can be seen from Table~\ref{tab:MeshQualityMonteCarlo} that the
probabilistically computed mesh obtained with $N=1000$ Monte Carlo simulations
yields relatively poor mesh quality measures when compared to the mesh obtained from the deterministic algorithm. If we increase the number of
Monte Carlo simulations to $N=10000$ the mesh quality measures found are already reasonably close to the deterministic case.

\subsection{Smoothing of the mesh}~\label{subsec:SmoothingOfMesh}

The problem with the stochastic representation of the solution of the linear mesh generator~\eqref{eq:WinslowTypeGridGeneratorStochasticSolution} is that usually a very large number of Monte Carlo simulations is needed to obtain a reasonably accurate numerical approximation. On the other hand, such a highly accurate numerical solution might not be absolutely necessary to obtain a good quality mesh from the grid generator~\eqref{eq:WinslowTypeGridGenerator}. Instead, we show that applying a smoothing operation
to a lower  accuracy solution of~\eqref{eq:WinslowTypeGridGenerator}
can yield a high-quality grid.

A type of smoothing algorithm that appears particularly well-suited for this type of problem is \textit{anisotropic diffusion}, also referred to as \emph{Perona--Malik diffusion},
\begin{align}\label{eq:PeronaMalikDiffusion}
\begin{split}
 &\xi_t-\nabla\cdot(c_\xi\nabla \xi)=0,\quad c_\xi=\exp(-||\nabla \xi||^2/k^2),\\
 &\eta_t-\nabla\cdot(c_\eta\nabla \eta)=0,\quad c_\eta=\exp(-||\nabla \eta||^2/k^2),
\end{split}
\end{align}
where $k$ is an arbitrary constant~\cite{pero90a}. In image processing, smoothing of this kind is used to remove noise without significantly affecting the edges of an image. This smoothing operation is suitable for grids generated with the Monte Carlo technique as the error in approximating the solution to~\eqref{eq:WinslowTypeGridGenerator} will appear as superimposed random noise over the required adapted mesh (the underlying signal). Anisotropic diffusion of the above form will remove this superimposed noise while still accurately preserving the regions of grid concentration.

\begin{figure}[!ht]
 \centering
 \includegraphics[scale=0.7]{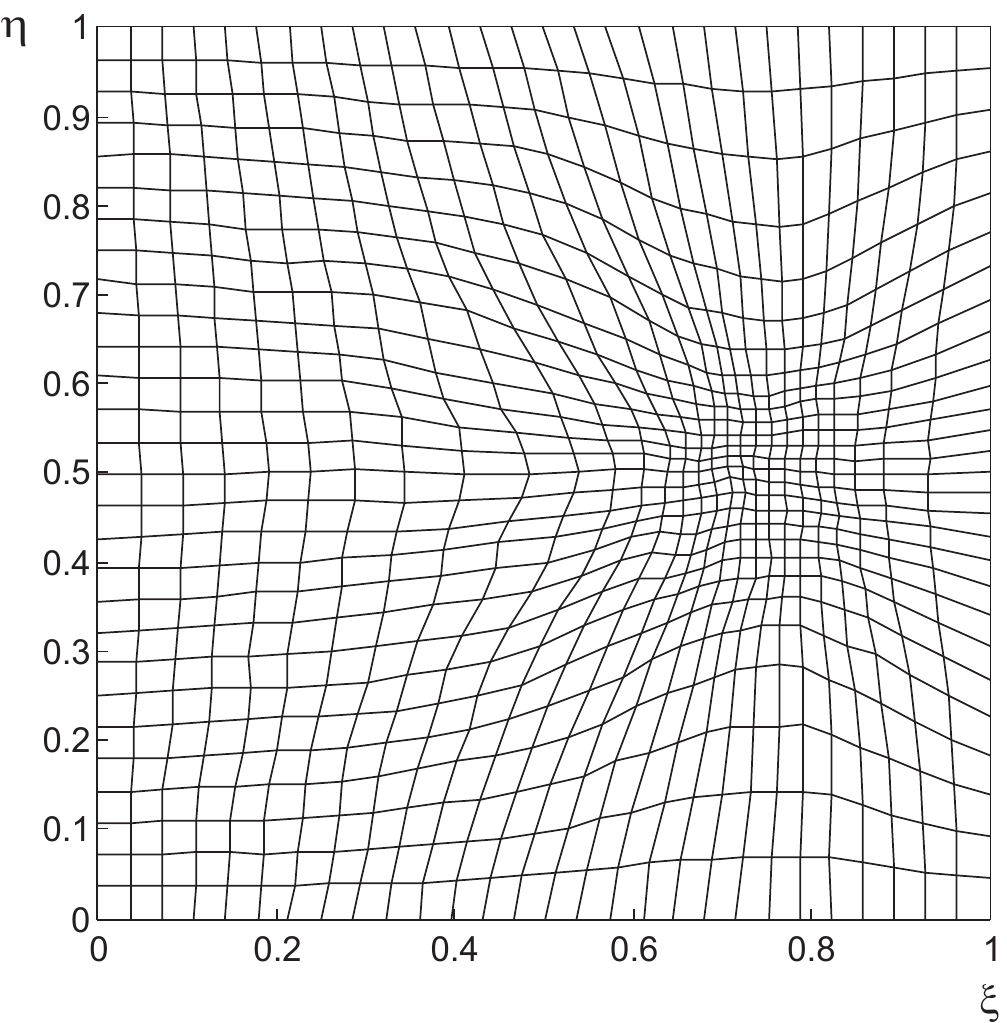}
 \caption{Grid obtained from the probabilistic solution~\eqref{eq:WinslowTypeGridGeneratorStochasticSolution} using $N=1000$ Monte Carlo simulations at each point to estimate the expected value. After the computation, we apply anisotropic diffusion of the form~\eqref{eq:PeronaMalikDiffusion} to the mesh in the computational space. Eq.~\eqref{eq:PeronaMalikDiffusion} was discretized using a FTCS scheme with $\Delta t=10^{-4}$ and $k=1000$. We integrate~\eqref{eq:PeronaMalikDiffusion} up to $t=5\Delta t$.}
 \label{fig:SingleDomainSolutionFullyProbabilistic10000SmoothedRunningExample}
\end{figure}

We now apply anisotropic diffusion to the fully probabilistically computed mesh displayed in Figure~\ref{fig:SingleDomainSolutionFullyProbabilisticRunningExample} (left). For this purpose we discretize~\eqref{eq:PeronaMalikDiffusion} with a forward in time and centered in space (FTCS) scheme. It is usually sufficient to numerically integrate~\eqref{eq:PeronaMalikDiffusion} only a few steps starting with a grid computed using the Monte Carlo technique to yield a smooth mesh. Consequently, diffusive smoothing is computationally a rather cheap operation.
The result of integrating~\eqref{eq:PeronaMalikDiffusion} to $t=5\Delta t$ starting with the
mesh shown on the left of
Figure~\ref{fig:SingleDomainSolutionFullyProbabilisticRunningExample} using a FTCS discretization with $k=1000$ and $\Delta t=10^{-4}$ is shown
in  Figure~\ref{fig:SingleDomainSolutionFullyProbabilistic10000SmoothedRunningExample}.

As can be seen by comparing with the single domain solution shown in the right of Figure~\ref{fig:SingleDomainSolutionRunningExample},  a few smoothing operations are able to almost completely eliminate the 'wiggles' that are typical of the approximate solution computed using Monte Carlo techniques. At the same time, the grid is still properly concentrating near the maximum values of the monitor function~\eqref{eq:MonitorFunctionRunningExample}, i.e., as expected, anisotropic diffusion properly preserves the signal in the Monte Carlo solution.

To confirm these qualitative findings we once again compute the geometric mesh quality measure for the probabilistically computed and smoothed mesh. Here it is particularly instructive to monitor the mesh quality as a function of the number of total smoothing steps $m$ (quantified through the final integration time $t=m\Delta t$). To have a proper comparison, we also smoothed the deterministic solution with the respective number of smoothing steps, which leads to the different absolute mesh quality measures for $Q_{\rm max}$ and $Q_{\rm mean}$ in Table~\ref{tab:MeshQualitySmoothedMesh}.

\begin{table}[!ht]
  \centering
  \caption{Mesh quality for grids adapting to monitor function~\eqref{eq:MonitorFunctionRunningExample}, smoothed using~\eqref{eq:PeronaMalikDiffusion} for $N=1000$ and $\lambda=1000$.}
  \begin{tabular}{|r|r|r|r|r|r|r|r|}
  \hline
   $m$ & $l_{\infty}$ & $Q_{\rm max}$ & $Q_{\rm mean}$ & $Q^{\rm S}_{\rm max}$ & $Q^{\rm S}_{\rm mean}$ & $R_{\rm max}$ & $R_{\rm mean}$ \\
   \hline
   $1$ & $2.25$ & $1.77$ & $1.16$ & $2.45$ & $1.2$ & $0.71$ & $0.97$ \\
   $2$ & $1.83$ & $1.77$ & $1.16$ & $2.44$ & $1.2$ & $0.72$ & $0.97$ \\
   $5$ & $1.72$ & $1.66$ & $1.15$ & $1.83$ & $1.18$ & $0.91$ & $0.98$ \\
   $10$ & $1.19$ & $1.54$ & $1.15$ & $1.57$ & $1.16$ & $0.98$ & $0.99$ \\
   \hline
   \end{tabular}
   \label{tab:MeshQualitySmoothedMesh}
\end{table}

\subsection{Domain decomposition solution}

Unless a large number of cores, equal to the total number of mesh points, is available,
the stochastic solution of (\ref{eq:WinslowTypeGridGenerator}) at all points would still be quite expensive.
In the context of domain decomposition we only evaluate the stochastic form of the solution at a few mesh points along the domain interfaces and then compute the solution deterministically at all the remaining grid points. As will be shown below, the meshes obtained from this stochastic--deterministic domain decomposition technique are usually much smoother than the meshes shown in Figure~\ref{fig:SingleDomainSolutionFullyProbabilisticRunningExample}. Thus, fewer (or no) smoothing sweeps will be required.

We now determine solution of the mesh generator~\eqref{eq:WinslowTypeGridGenerator} with the monitor function~\eqref{eq:MonitorFunctionRunningExample} using the domain decomposition technique outlined in Section~\ref{sec:TheoryDD}.  Here we restrict ourselves to the case of four square subdomains. The probabilistic solutions at the interfaces are computed using $\lambda=10000$ and $N=10000$.
Additional experiments with larger number of subdomains are given in the next section.

\begin{figure}[!ht]
 \centering
 \includegraphics[scale=0.7]{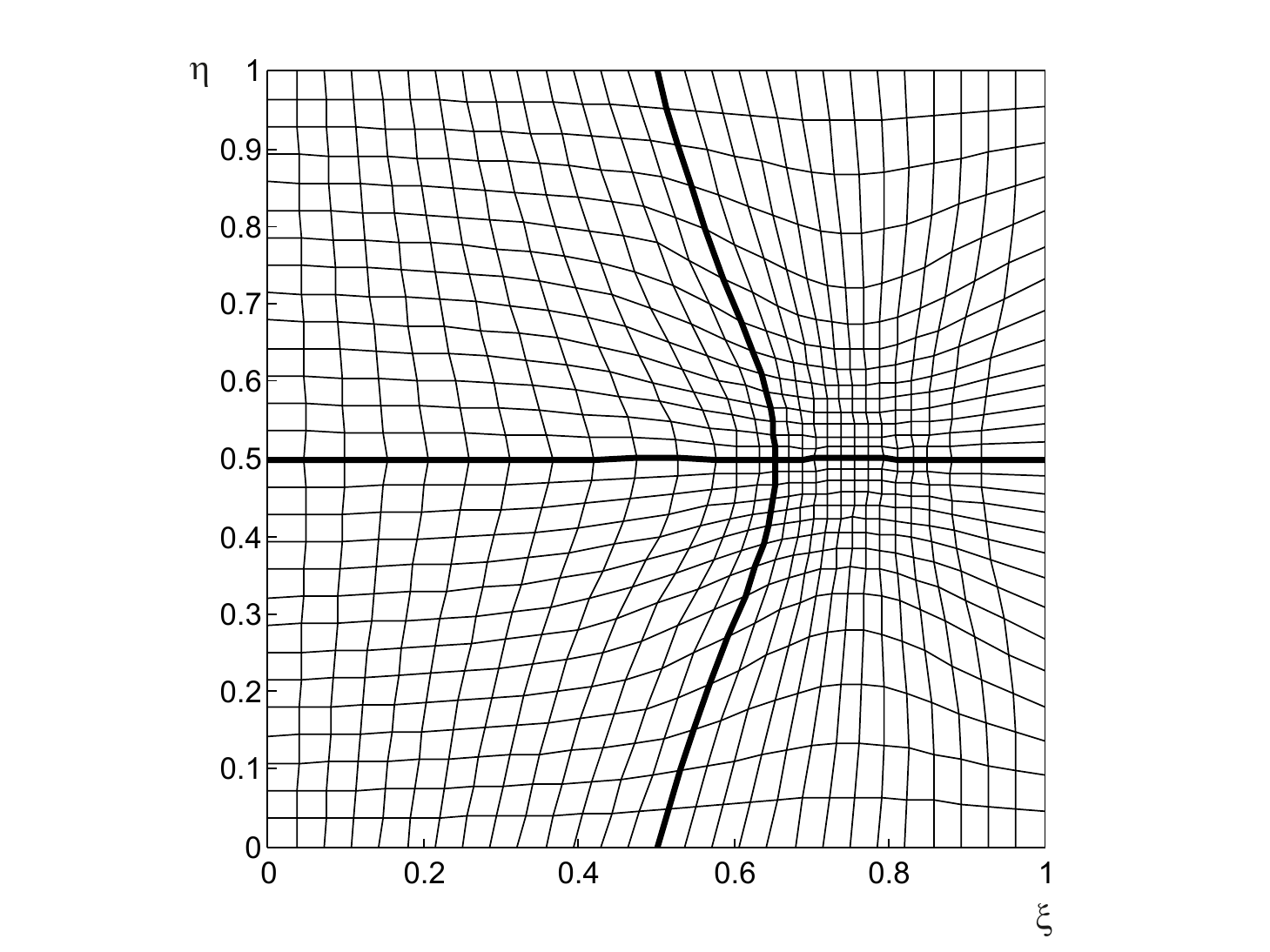}
 \caption{Domain decomposition solution for the mesh generator~\eqref{eq:WinslowTypeGridGenerator} with monitor function~\eqref{eq:MonitorFunctionRunningExample}.}
 \label{fig:DDSolutionRunningExample}
\end{figure}

As can be seen from Figure~\ref{fig:DDSolutionRunningExample} the stochastic--deterministic approach to solve~\eqref{eq:WinslowTypeGridGenerator} almost perfectly reproduces the single domain solution shown in Figure~\ref{fig:SingleDomainSolutionRunningExample} (right).
 This finding is also confirmed with the geometric mesh quality measure, having values $Q^{\rm DD}_{\rm max}=1.80$ and $Q^{\rm DD}_{\rm mean}=1.16$ for the domain decomposition solution which give rise to the ratios $R_{\rm max}=1$ and $R_{\rm mean}=1$ when compared to the single domain solution.

No smoothing was applied to obtain the mesh shown in Figure ~\ref{fig:DDSolutionRunningExample}. We should nevertheless like to stress that the smoothing operation proposed in Section~\ref{subsec:SmoothingOfMesh}, if needed,  can be also achieved in parallel, as~\eqref{eq:PeronaMalikDiffusion} can be applied on each subdomain separately. That is, equation~\eqref{eq:PeronaMalikDiffusion} yields a local smoothing operation and does not dilute the parallel nature of the proposed algorithm.

The only complication is how to smooth the subdomain interfaces as they serve as natural (Dirichlet) boundaries for the smoothing operation~\eqref{eq:PeronaMalikDiffusion}. Ways to overcome this difficulty include a pre-smoothing of the interface to ensure that it is sufficiently accurate before being used as boundary conditions for the subdomain solves
or possibly apply a second smoothing cycle with shifted subdomains that have the original interfaces in their interior. We will present an example for the first possibility in Section~\ref{sec:ResultsDD}.

\subsection{Interpolation along the interface}

If the number of compute cores is limited, a  promising approach to further reduce the cost of the probabilistic part of the domain decomposition algorithm is to avoid computing the solution stochastically at all the grid points along the interface.
Instead we obtain the stochastic solution at only a few points and approximate the solution at the remaining points using interpolation.  This was suggested in \cite{aceb05a}, however,
we found the meshes obtained by this approach to be quite
sensitive to the location of these points.

\begin{figure}[!ht]
\centering
\begin{subfigure}[b]{0.45\textwidth}
  \centering
  \includegraphics[width=\linewidth]{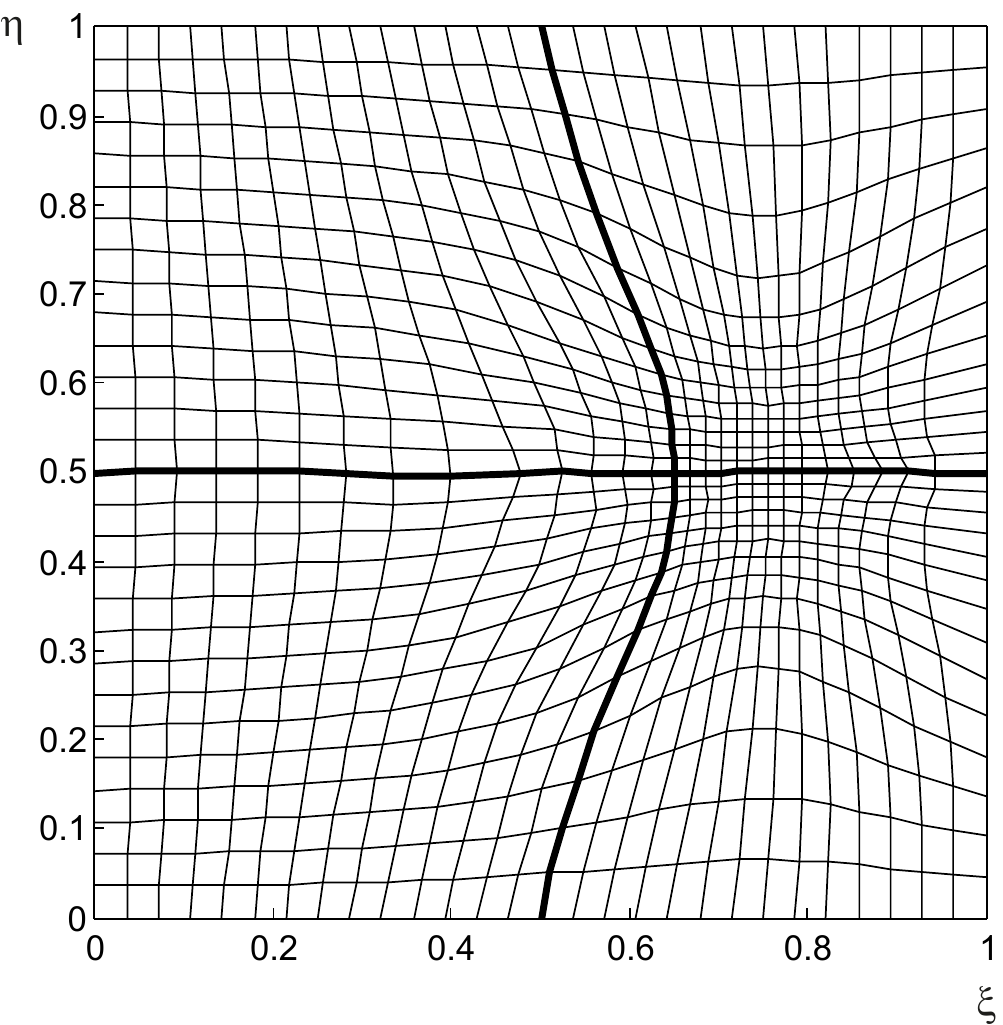}
  \label{fig:DDSolutionNoPlacementRunningExample}
\end{subfigure}\qquad
\begin{subfigure}[b]{0.45\textwidth}
  \centering
  \includegraphics[width=\linewidth]{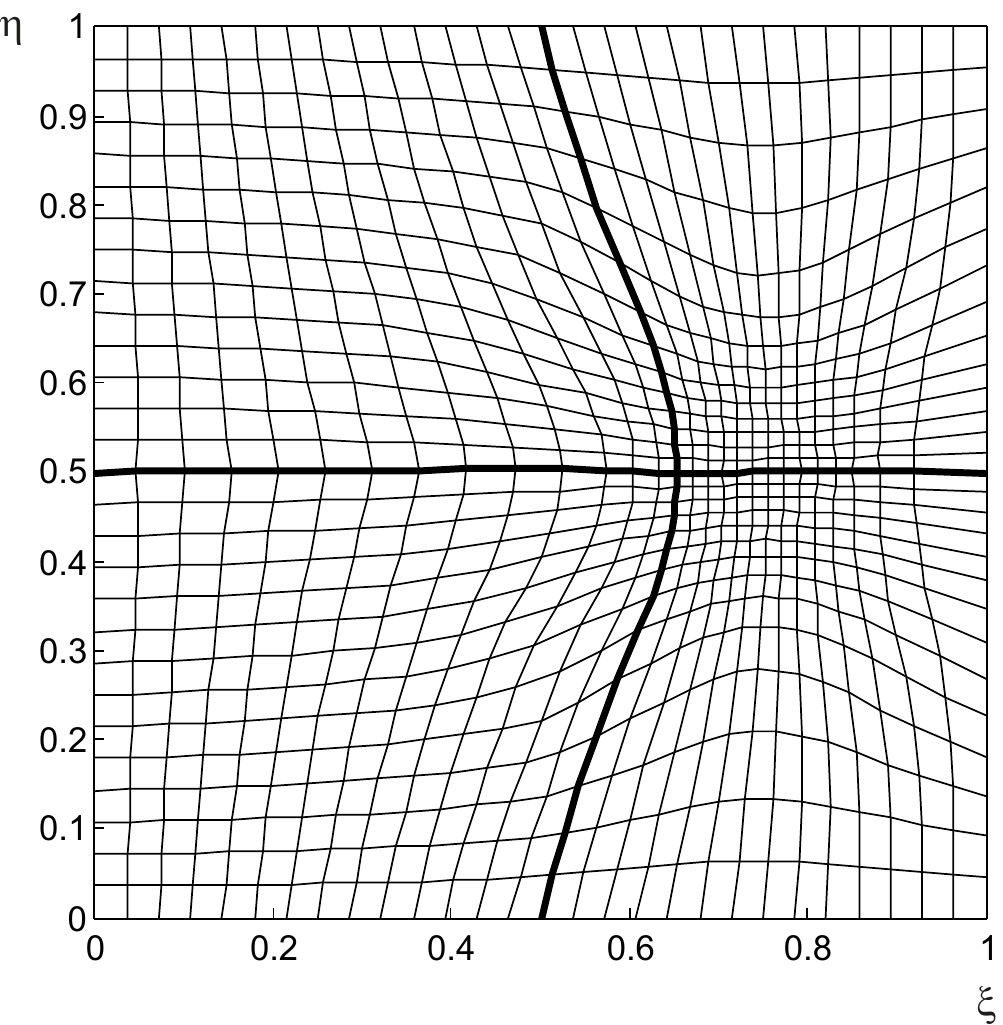}
  \label{fig:DDSolutionWithPlacementRunningExample}
\end{subfigure}
\caption{\textbf{Left:} Domain decomposition solution for the mesh generator~\eqref{eq:WinslowTypeGridGenerator} with monitor function~\eqref{eq:MonitorFunctionRunningExample}. The solution was evaluated at seven equidistant points along each dividing line probabilistically. Linear interpolation was used to obtain the remaining interface points. \textbf{Right:} Same as left, but solution was evaluated probabilistically at the maxima and minima of the first and second derivatives of the monitor function at the interface. This gave seven points where the mesh was computed with Monte Carlo simulations. Linear interpolation was used to obtain the remaining interface points.}
\label{fig:DDSolutionInterpolationInterfaceRunningExample}
\end{figure}

A short study in this \emph{optimal placement problem} is reported in Figure~\ref{fig:DDSolutionInterpolationInterfaceRunningExample}. In Figure~\ref{fig:DDSolutionInterpolationInterfaceRunningExample} (left) we computed the probabilistic solution at seven equally spaced points along each of the dividing lines. Linear interpolation is used to obtain the remaining interface points. It is clearly visible that this procedure does not give the mesh obtained by domain decomposition in Figure~\ref{fig:DDSolutionRunningExample}. The problem is that local maxima and minima in the monitor function is
missed and hence proper transitions between regions of grid concentration and de-concentration are not captured. We also tested more sophisticated interpolation strategies such as splines, cubic interpolation and Chebychev interpolation, but the results obtained are nearly identical.
We conclude that, unlike the results reported in \cite{aceb05a}, the placement of points is crucial if interpolation is used to approximate the interface conditions for the domain decomposition for mesh generation.

In Figure~\ref{fig:DDSolutionInterpolationInterfaceRunningExample} (right) we determined the points along the  artificial interfaces where the solution is computed probabilistically based on the properties of the monitor function $\rho$. In particular, we required placement of the points near the maxima and minima of $\rho_x$ and $\rho_{xx}$ along the horizontal dividing line and $\rho_y$ and $\rho_{yy}$ along the vertical dividing line. This leads to a total of seven points along each  dividing line where the solution has to be computed using the Monte Carlo technique. Note that these points are not equally spaced as they were in the previous case. The remaining interface points are then approximated using linear interpolation.

The solution found using this rather simple placement criterion gives a mesh very close to the solution shown in Figure~\ref{fig:DDSolutionRunningExample} but at only a fraction of the computational cost, or cores, required by the original stochastic domain decomposition solution. Smoothing of the mesh shown could be used to further improve the smoothness in the transition of the grid lines across the interfaces.

These qualitatively visible differences in Figure~\ref{fig:DDSolutionInterpolationInterfaceRunningExample} are not clearly reflected in our chosen mesh quality measure, which yields essentially the same results ($R=1$) for both equally and 'optimally' placed interfaced points.
On the other hand, the $l_\infty$ error found when using equally spaced points, $l_\infty=4.91$, is more than five times as large as the $l_\infty$ error
obtained when we strategically place the points where the stochastic solution
is obtained, $l_\infty=0.92$.  This  shows a larger deviation from the single domain solution.

The one drawback of this placement approach is that it is not obvious a priori how to determine how many and at what locations the solution is best determined probabilistically. For large meshes in real-world applications, a possible trade-off would be to place the points near the most pronounced features of the monitor function only and smooth the resulting meshes over a few cycles.

\section{Further examples}\label{sec:ResultsDD}

Having explained the technique and the issues involved, in this section we present additional examples to demonstrate the meshes generated using our stochastic domain decomposition algorithm on well known test problems.

We first choose the monitor function
\begin{equation}\label{eq:HuangSloanMonitorFunction}
 \rho=\sqrt {1+\alpha(R^2\exp(R(x-1))^2\sin(\pi y)^2+(1-\exp(R(x-1)))^2\cos(\pi y)^2\pi^2 )}
\end{equation}
with $\alpha=0.7$ and $R=15$.   This monitor function was used in~\cite{hayn12a,huan94a} to generate meshes using a nonlinear mesh generator.

\begin{figure}[!ht]
\centering
\begin{subfigure}[b]{0.45\textwidth}
  \centering
  \includegraphics[width=\linewidth]{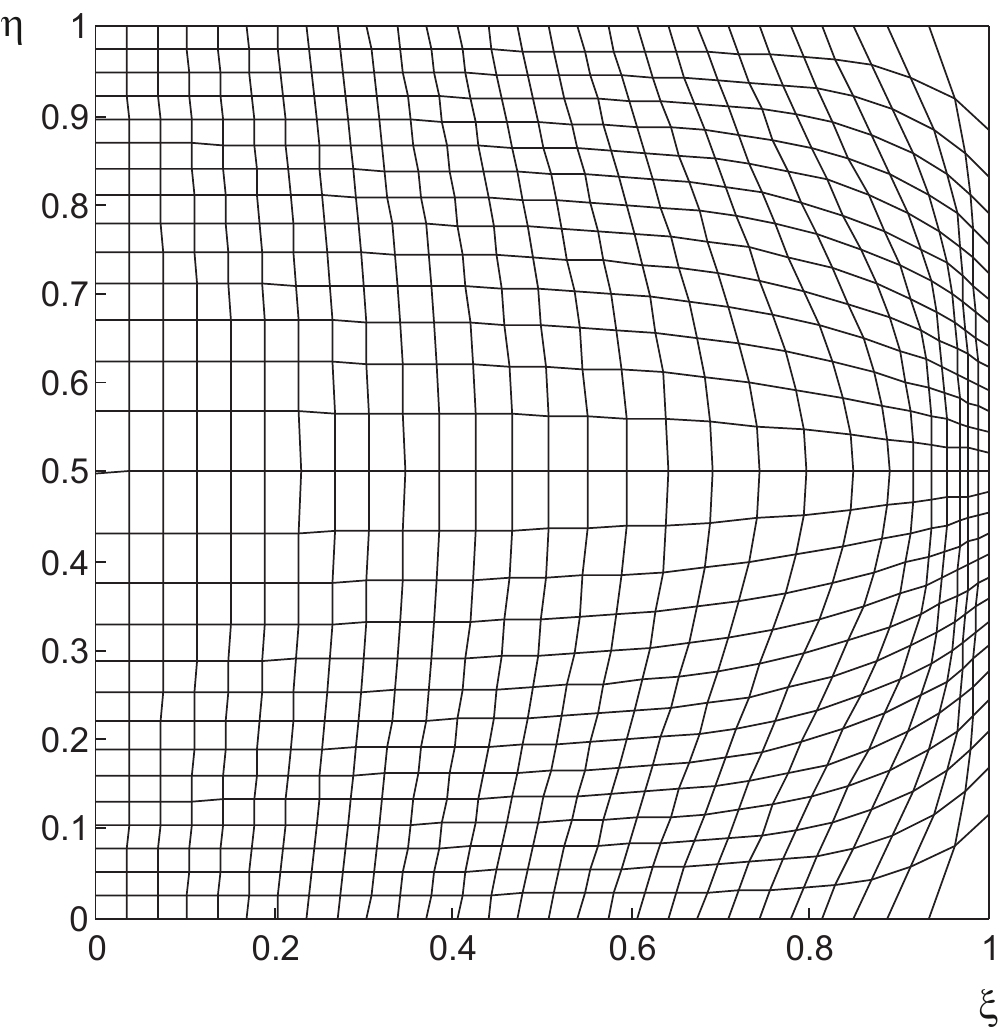}
  \label{fig:HuangSloanSingleDomain}
\end{subfigure}\qquad
\begin{subfigure}[b]{0.45\textwidth}
  \centering
  \includegraphics[width=\linewidth]{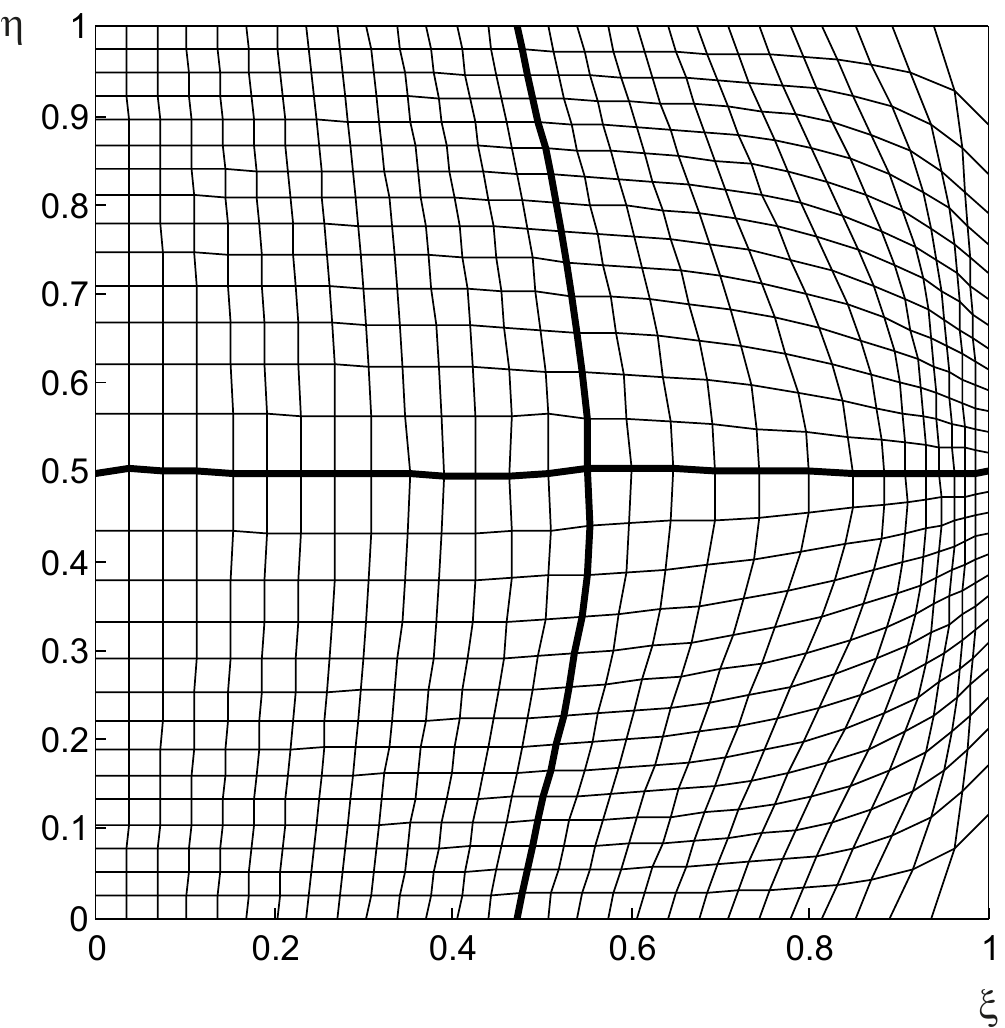}
  \label{fig:HuangSloanDD}
\end{subfigure}
\caption{\textbf{Left:} Single domain solution of the mesh generator~\eqref{eq:WinslowTypeGridGenerator} using the monitor function~\eqref{eq:HuangSloanMonitorFunction}. \textbf{Right:} The associated domain decomposition solution. The solution along the interfaces was computed using~\eqref{eq:WinslowTypeGridGeneratorStochasticSolution} with $N=10000$ Monte Carlo simulations and $\lambda=10000$ as the parameter of the exponential distribution.}
\label{fig:HuangSloan}
\end{figure}

Figure~\ref{fig:HuangSloan} (left) displays the single domain mesh obtained by solving
Winslow's generator~\eqref{eq:WinslowTypeGridGenerator} using a Jacobi iteration. The corresponding domain decomposition solution with four subdomains is depicted in Figure~\ref{fig:HuangSloan} (right). Here, we used $\lambda=10000$ and $N=10000$ to compute the solution probabilistically along all the points on the interfaces. No smoothing was applied to the final mesh.

 The mesh quality measures for the domain decomposition solution are $Q^{\rm DD}_{\rm max}=1.69$ and $Q^{\rm DD}_{\rm mean}=1.14$, respectively, with corresponding ratios $R_{\rm max}=0.99$ and $R_{\rm mean}=1$ compared to the single domain solution. That is, in terms of the geometric mesh quality measure, the single domain solution is almost perfectly approximated.

\bigskip

\noindent The next example uses the monitor function
\begin{equation}\label{eq:MackenzieMonitorFunction}
 \rho=\frac{1}{1+\alpha\exp(-R(y-1/2-1/4\sin(2\pi x))^2)}
\end{equation}
with $\alpha=10$ and $R=50$.

\begin{figure}[!ht]
\centering
\begin{subfigure}[b]{0.45\textwidth}
  \centering
  \includegraphics[width=\linewidth]{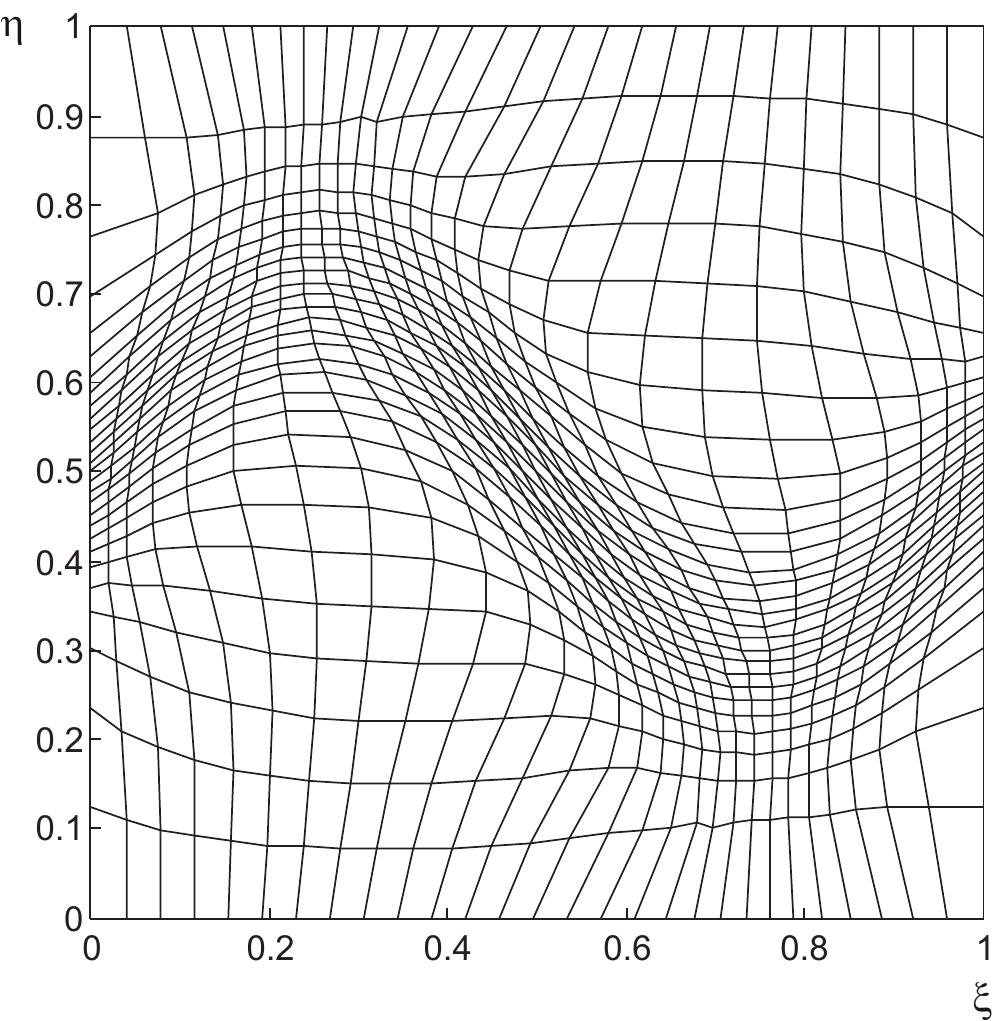}
  \label{fig:MackenzieSingleDomain}
\end{subfigure}\qquad
\begin{subfigure}[b]{0.45\textwidth}
  \centering
  \includegraphics[width=\linewidth]{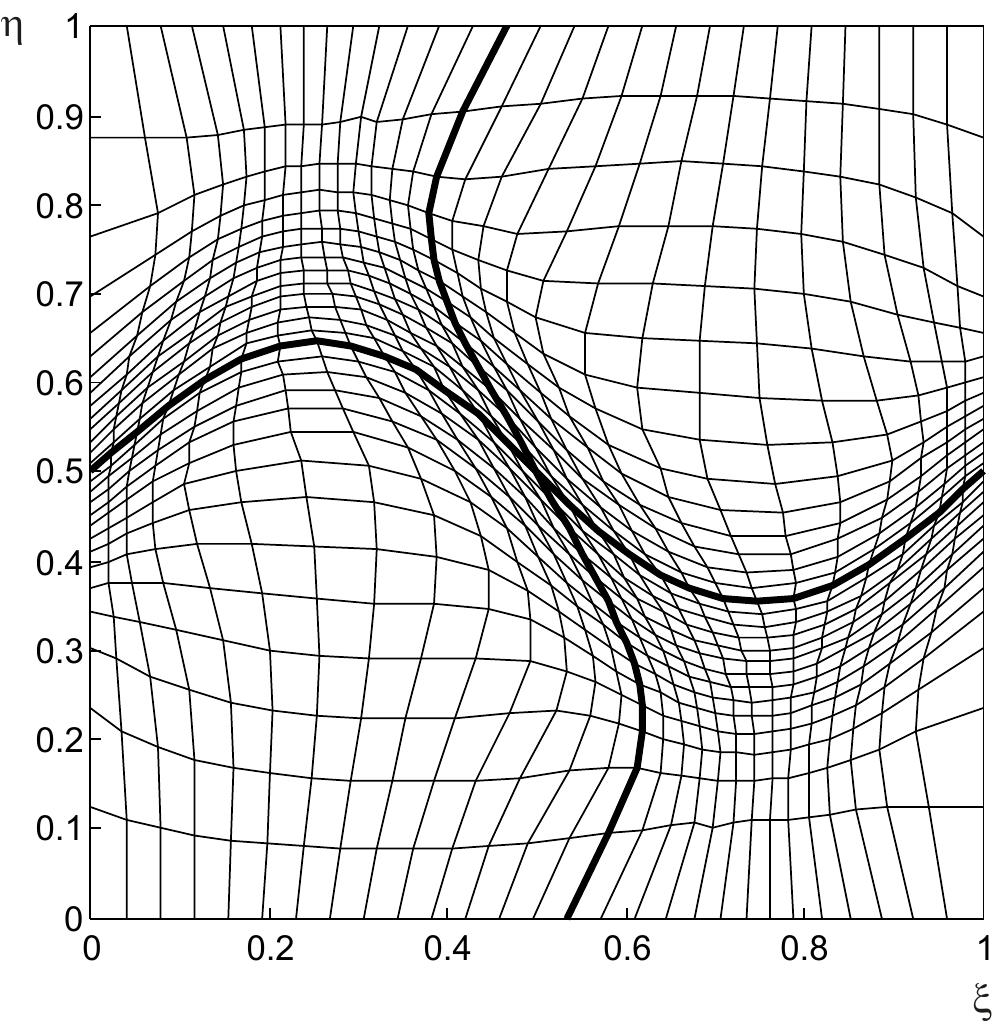}
  \label{fig:MackenzieDD}
\end{subfigure}
\caption{\textbf{Left:} Single domain solution of the mesh generator~\eqref{eq:WinslowTypeGridGenerator} using the monitor function~\eqref{eq:MackenzieMonitorFunction}. \textbf{Right:} The associated domain decomposition solution. The solutions at the interface points were computed using~\eqref{eq:WinslowTypeGridGeneratorStochasticSolution} with $N=10000$ Monte Carlo simulations and $\lambda=10000$ as the parameter of the exponential distribution.}
\label{fig:Mackenzie}
\end{figure}

The single domain solution for this example is depicted in Figure~\ref{fig:Mackenzie} (left).
The mesh reflects both
the large scale features (the sinusoidal wave) and smaller scale variations (across the wave) in the
monitor function.
In Figure~\ref{fig:Mackenzie} (right) we show the domain decomposition solution using four subdomains. The parameters for the stochastic solver again were $\lambda=10000$ and $N=10000$. We used the stochastic solver to compute all the points along the interfaces. No smoothing is applied to the final mesh.
It can be seen in Figure~\ref{fig:Mackenzie} (right) that the domain decomposition solution has several kinks which are slightly more pronounced than in the single
domain mesh.

\begin{figure}[!ht]
\centering
\begin{subfigure}[b]{0.45\textwidth}
  \centering
  \includegraphics[width=\linewidth]{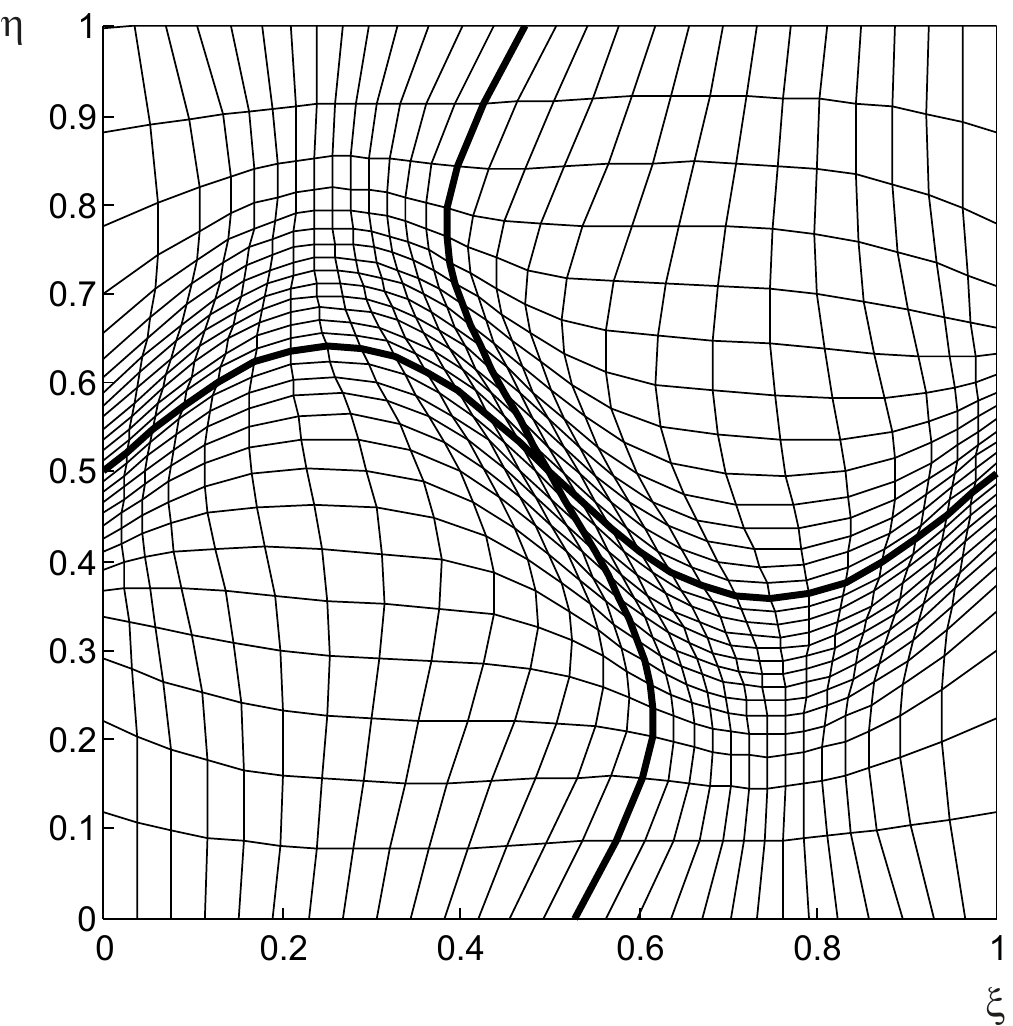}
\end{subfigure}\qquad
\begin{subfigure}[b]{0.45\textwidth}
  \centering
  \includegraphics[width=\linewidth]{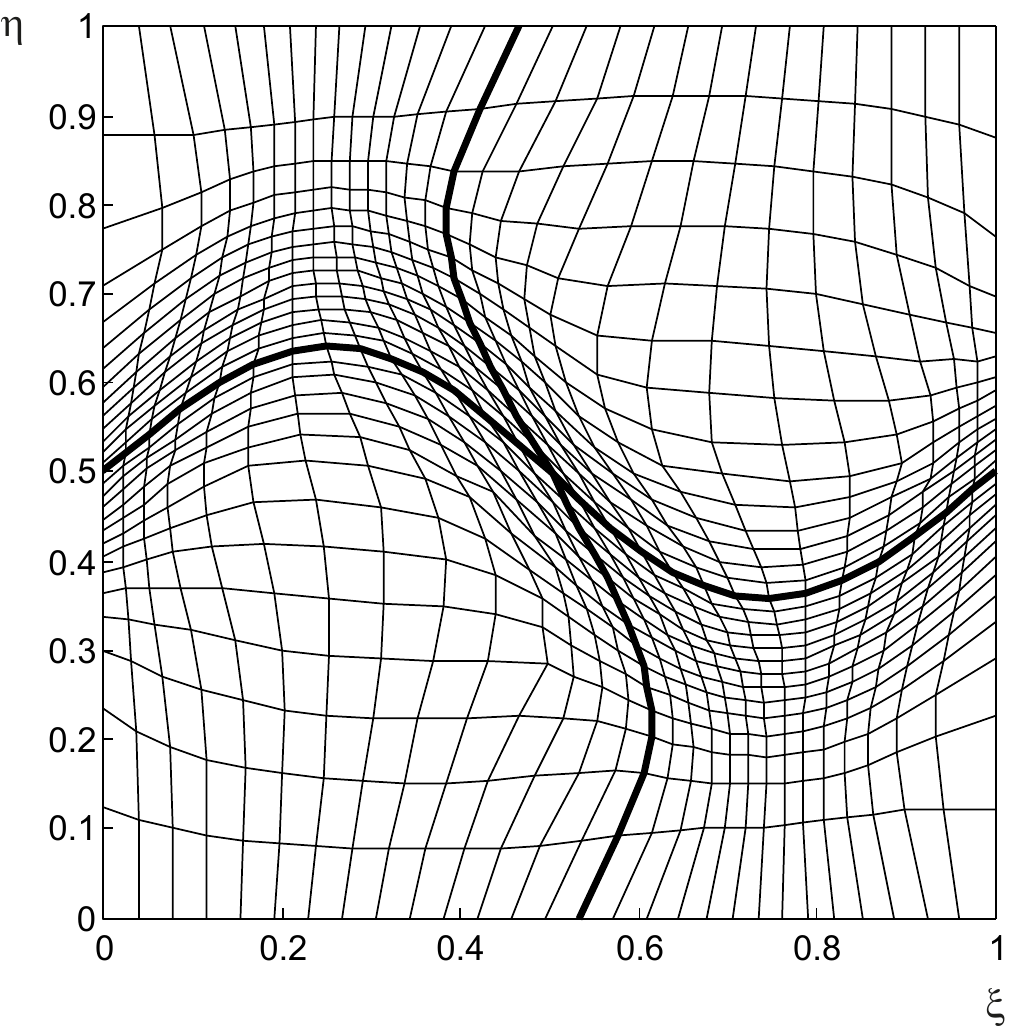}
\end{subfigure}
\caption{\textbf{Left:} Domain decomposition solution of the mesh generator~\eqref{eq:WinslowTypeGridGenerator} using the monitor function~\eqref{eq:MackenzieMonitorFunction}. Five smoothing cycles using~\eqref{eq:PeronaMalikDiffusion} with $\Delta t=10^{-4}$ and $k=1000$ were applied on the whole domain. \textbf{Right:} Smoothing is applied on each of the subdomains separately.}
\label{fig:MackenzieSmoothedGlobalSubdomains}
\end{figure}

\begin{figure}[!ht]
\centering
\begin{subfigure}[b]{0.45\textwidth}
  \centering
  \includegraphics[width=\linewidth]{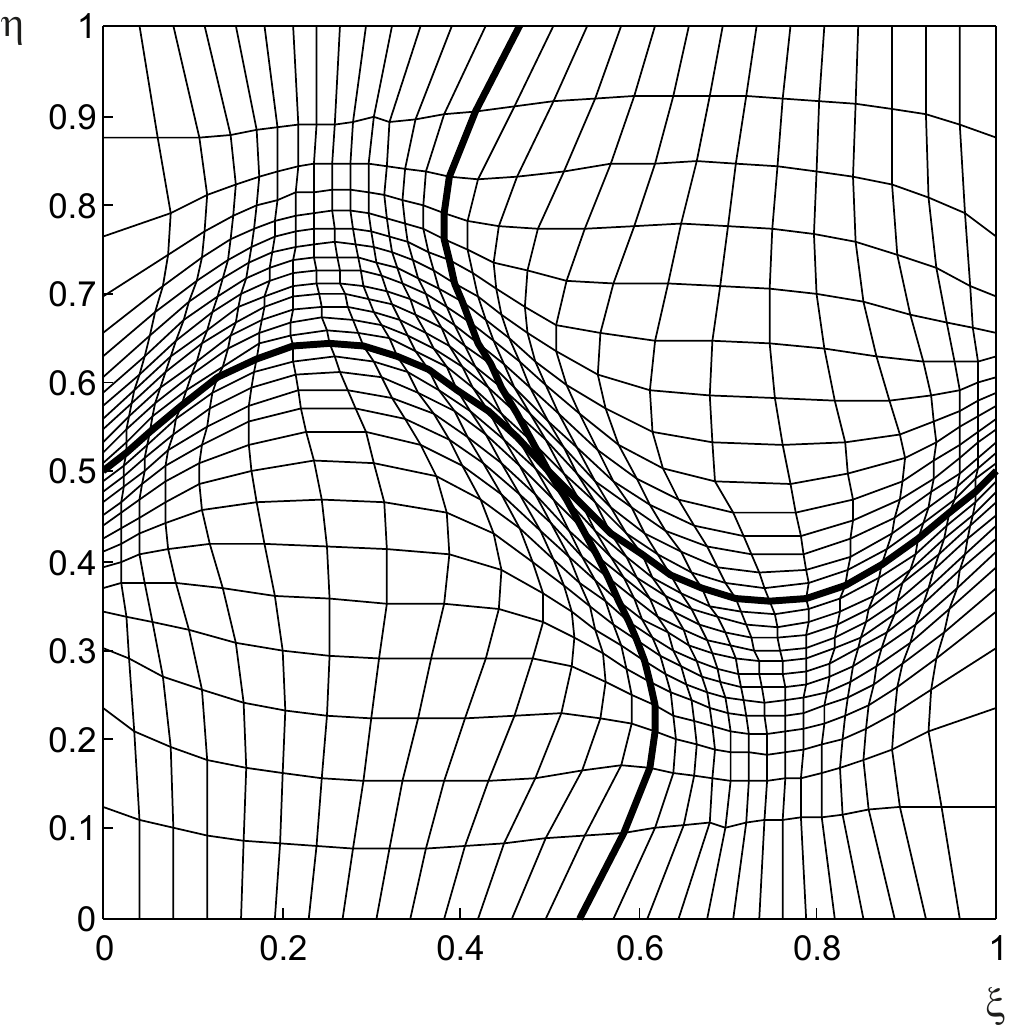}
  \label{fig:MackenzieSmoothedInterface}
\end{subfigure}\qquad
\begin{subfigure}[b]{0.45\textwidth}
  \centering
  \includegraphics[width=\linewidth]{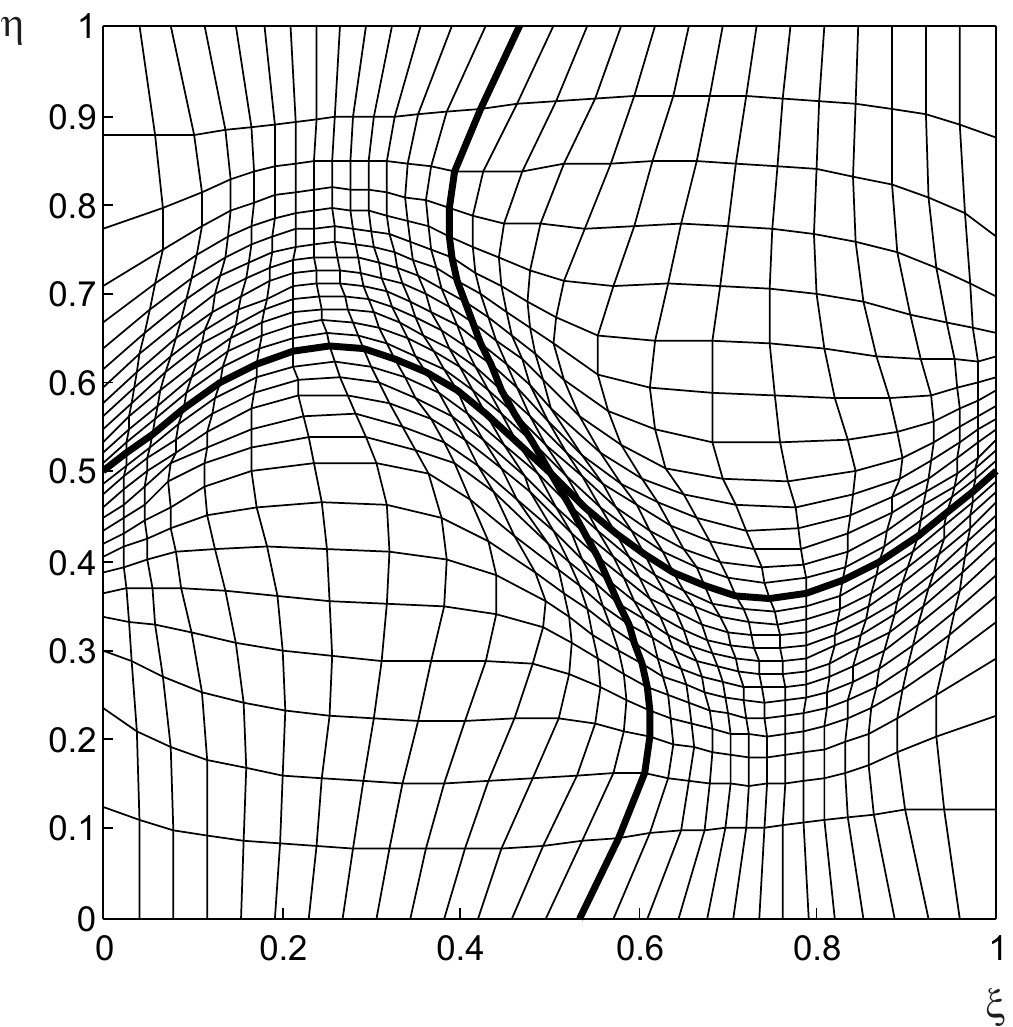}
  \label{fig:MackenzieSmoothedInterfaceSubdomains}
\end{subfigure}
\caption{\textbf{Left:} Domain decomposition solution of the mesh generator~\eqref{eq:WinslowTypeGridGenerator} using the monitor function~\eqref{eq:MackenzieMonitorFunction}. The probabilistically computed interface solution was smoothed using weighted linear least squares and a second degree polynomial model. \textbf{Right:} Same as right, but additionally five smoothing cycles using~\eqref{eq:PeronaMalikDiffusion} with $\Delta t=10^{-4}$ and $k=1000$ were applied on each of the subdomains separately.}
\label{fig:MackenzieDDSmoothedInterface}
\end{figure}

We now present the results of four different versions of the domain decomposition solution in Figure~\ref{fig:MackenzieSmoothedGlobalSubdomains} and Figure~\ref{fig:MackenzieDDSmoothedInterface}.
In Figure~\ref{fig:MackenzieSmoothedGlobalSubdomains} we contrast the effect of global smoothing (left) versus local smoothing on each of the subdomains (right). As was indicated in Section~\ref{sec:DDStochastic} the global smoothing is a sequential operation whereas smoothing on the subdomains can be carried out in parallel.  The global smoothing operation leads to a mesh that varies very smoothly throughout the whole domain.   The locally smoothed mesh is much improved as well (also compare the mesh quality measures reported in Tables~\ref{tab:MeshQualityMeasuresConstantN} and \ref{tab:MeshQualityMeasuresConstantLambda}) and has
less kinks when compared to the original domain decomposition solution shown in Figure~\ref{fig:Mackenzie} (right). The remaining discontinuities could be further reduced by applying a second smoothing cycle over shifted subdomains, i.e.\ by re-assigning the subdomains in such a manner that the original interfaces lie within the new subdomains and then re-apply the smoothing. Once again, this would not sacrifice the overall parallel nature of the proposed algorithm.

In Figure~\ref{fig:MackenzieDDSmoothedInterface} we demonstrate the effect of pre-smoothing the stochastically computed interface solution.  In Figure~\ref{fig:MackenzieDDSmoothedInterface} (left), the
values of $\xi$ and $\eta$ have been smoothed using weighted linear least squares and a second degree polynomial model (provided by the
\texttt{Matlab} function \texttt{smooth}) before using the interface solution as the boundary values for the solver on each of the subdomains.
In general  we see that only a few interface solutions are required; smoothing of the interface is much cheaper than smoothing the whole domain or the subdomains. In Figure~\ref{fig:MackenzieDDSmoothedInterface} (right), we demonstrate the effect of both smoothing the interface and applying the diffusive smoothing~\eqref{eq:PeronaMalikDiffusion} on each of the subdomains.

 We now present the results for the geometric mesh quality measure for the different variants of the stochastic--deterministic domain decomposition algorithm discussed above. In particular, the geometric mesh quality measures for the grids computed using monitor function~\eqref{eq:MackenzieMonitorFunction} are determined for two series of numerical experiments. In Table~\ref{tab:MeshQualityMeasuresConstantN}, we keep the number of Monte Carlo runs constant, $N=10000$, and vary the parameter of the exponential distribution~$\lambda$; in  Table~\ref{tab:MeshQualityMeasuresConstantLambda} we vary~$N$ and keep $\lambda=5000$ constant.


\begin{table}[!ht]
  \centering
  \caption{Mesh quality for grids adapting to monitor function~\eqref{eq:MackenzieMonitorFunction} with $N=10000$ Monte Carlo simulations.}
  \begin{tabular}{|p{4cm}|r|r|r|r|r|r|}
  \hline
   & $\lambda$ & $l_{\infty}$ & $Q^{\rm DD}_{\rm max}$ & $Q^{\rm DD}_{\rm mean}$ & $R_{\rm max}$ & $R_{\rm mean}$ \\
   \hline
   \multirow{4}{*}{\parbox{4cm}{No smoothing\\ $Q_{\rm max}=3.86$\\ $Q_{\rm mean}=1.67$}}& $1000$ & $1.18$ & $4.44$ & $1.76$ & $0.87$ & $0.95$ \\
   &$2000$ & $0.91$ & $4.29$ & $1.72$ & $0.90$ & $0.97$ \\
   &$5000$ & $0.49$ & $3.98$ & $1.7$ & $0.97$ & $0.98$ \\
   &$10000$ & $0.4$ & $3.9$ &  $1.69$ & $0.99$ & $0.99$ \\
  \hline
  \multirow{4}{*}{\parbox{4cm}{Global smoothing\\ $Q_{\rm max}=3.6$\\ $Q_{\rm mean}=1.62$}}& $1000$ & $1.04$ & $3.96$ & $1.71$ & $0.91$ & $0.95$ \\
   &$2000$ & $0.64$ & $3.79$ & $1.67$ & $0.95$ & $0.97$ \\
   &$5000$ & $0.45$ & $3.67$ & $1.64$ & $0.98$ & $0.99$ \\
   &$10000$ & $0.44$ & $3.67$ & $1.62$ & $0.98$ & $1$ \\
  \hline
  \multirow{4}{*}{\parbox{4cm}{Sub-domain smoothing\\ $Q_{\rm max}=3.6$\\ $Q_{\rm mean}=1.62$}} & $1000$ & $1.36$ & $4.24$ & $1.69$ & $0.85$ & $0.96$ \\
   &$2000$ & $0.92$ & $3.83$ & $1.65$ & $0.94$ & $0.98$ \\
   &$5000$ & $0.75$ & $3.79$ & $1.62$ & $0.95$ & $1$ \\
   &$10000$ & $0.74$ & $3.67$ & $1.62$ & $0.98$ & $1$ \\
  \hline
  \multirow{4}{*}{\parbox{4cm}{Interface smoothing\\ $Q_{\rm max}=3.86$\\ $Q_{\rm mean}=1.67$}}& $1000$ & $1.40$ & $4.89$ & $1.76$ & $0.79$ & $0.95$ \\
   &$2000$ & $0.82$ & $4.54$ & $1.73$ & $0.85$ & $0.97$ \\
   &$5000$ & $0.5$ & $4.39$ & $1.69$ & $0.88$ & $0.99$ \\
   &$10000$ & $0.41$ & $4.29$ & $1.69$ & $0.9$ & $0.99$ \\
  \hline
  \multirow{4}{*}{\parbox{4cm}{Sub-domain and\\ interface smoothing\\ $Q_{\rm max}=3.6$\\ $Q_{\rm mean}=1.62$}} &$1000$ & $1.52$ & $4.29$ & $1.69$ & $0.84$ & $0.96$ \\
   &$2000$ & $1.1$ & $4.8$ & $1.65$ & $0.75$ & $0.98$ \\
   &$5000$ & $0.8$ & $4.04$ & $1.62$ & $0.89$ & $1$ \\
   &$10000$ & $0.75$ & $4.09$ & $1.62$ & $0.88$ & $1$ \\
  \hline
  \end{tabular}
  \label{tab:MeshQualityMeasuresConstantN}
\end{table}


\begin{table}[!ht]
  \centering
  \caption{Mesh quality for grids adapting to monitor function~\eqref{eq:MackenzieMonitorFunction} with  $\lambda=5000$.}
  \begin{tabular}{|p{4cm}|r|r|r|r|r|r|}
  \hline
   & $N$ & $l_{\infty}$ & $Q^{\rm DD}_{\rm max}$ & $Q^{\rm DD}_{\rm mean}$ & $R_{\rm max}$ & $R_{\rm mean}$ \\
   \hline
   \multirow{4}{*}{\parbox{4cm}{No smoothing\\ $Q_{\rm max}=3.86$\\ $Q_{\rm mean}=1.67$}}& $1000$ & $0.74$ & $5.44$ & $1.72$ & $0.71$ & $0.97$ \\
   &$2000$ & $0.66$ & $3.99$ & $1.69$ & $0.97$ & $0.99$ \\
   &$5000$ & $0.4$ & $3.9$ & $1.69$ & $0.99$ & $0.99$ \\
   &$10000$ & $0.45$ & $3.9$ & $1.72$ & $0.99$ & $0.97$ \\
  \hline
  \multirow{4}{*}{\parbox{4cm}{Global smoothing\\ $Q_{\rm max}=3.6$\\ $Q_{\rm mean}=1.62$}}&$1000$ & $0.762$ & $3.79$ & $1.65$ & $0.95$ & $0.98$ \\
   &$2000$ & $0.47$ & $3.83$ & $1.64$ & $0.94$ & $0.99$ \\
   &$5000$ & $0.36$ & $3.75$ & $1.64$ & $0.96$ & $0.99$ \\
   &$10000$ & $0.44$ & $3.67$ & $1.62$ & $0.98$ & $1$ \\
  \hline
  \multirow{4}{*}{\parbox{4cm}{Sub-domain smoothing\\ $Q_{\rm max}=3.6$\\ $Q_{\rm mean}=1.62$}}&$1000$ & $1.15$ & $4.5$ & $1.64$ & $0.8$ & $0.99$ \\
   &$2000$ & $0.94$ & $4$ & $1.64$ & $0.89$ & $0.99$ \\
   &$5000$ & $0.82$ & $3.96$ & $1.62$ & $0.91$ & $1$ \\
   &$10000$ & $0.74$ & $3.67$ & $1.62$ & $0.98$ & $1$ \\
  \hline
  \multirow{4}{*}{\parbox{4cm}{Interface smoothing\\ $Q_{\rm max}=3.86$\\ $Q_{\rm mean}=1.67$}}&$1000$ & $1.47$ & $4.82$ & $1.7$ & $0.8$ & $0.98$ \\
   &$2000$ & $0.95$ & $4.95$ & $1.7$ & $0.78$ & $0.98$ \\
   &$5000$ & $0.59$ & $4.95$ & $1.7$ & $0.78$ & $0.98$ \\
   &$10000$ & $0.41$ & $4.29$ & $1.69$ & $0.9$ & $0.99$ \\
  \hline
  \multirow{4}{*}{\parbox{4cm}{Sub-domain and\\ interface smoothing\\ $Q_{\rm max}=3.6$\\ $Q_{\rm mean}=1.62$}}&$1000$ & $1.1$ & $6$ & $1.65$ & $0.6$ & $0.98$ \\
   &$2000$ & $0.95$ & $4.86$ & $1.64$ & $0.74$ & $0.99$ \\
   &$5000$ & $0.82$ & $4.5$ & $1.64$ & $0.8$ & $0.99$ \\
   &$10000$ & $0.75$ & $4.01$ & $1.62$ & $0.88$ & $1$ \\
  \hline
  \end{tabular}
  \label{tab:MeshQualityMeasuresConstantLambda}
\end{table}

Tables~\ref{tab:MeshQualityMeasuresConstantN} and~\ref{tab:MeshQualityMeasuresConstantLambda} illustrate that both the value of $\lambda$, the exponential
time-stepping parameter, and the number of Monte Carlo simulations $N$ affect the quality of the resulting mesh.  The larger $\lambda$ and $N$, the better the domain decomposition mesh approximates the quality of the single domain mesh. This is expected since as
$\lambda\to\infty$ and $N\to\infty$ the numerical approximation of the stochastic representation of the solution to the linear mesh generator~\eqref{eq:WinslowTypeGridGenerator} approaches the deterministic solution and thus converges to the single domain result.

From Tables~\ref{tab:MeshQualityMeasuresConstantN} and~\ref{tab:MeshQualityMeasuresConstantLambda} it can be seen that the meshes obtained using domain decomposition with probabilistically computed solution along the interfaces give  good quality meshes (compared to the single domain solution) with very few Monte Carlo simulations (low $N$) at relatively coarse mean time steps (low $\lambda$). Global and local smoothing on the subdomains can bring the mesh quality of the domain decomposition solution slightly closer to the single domain case, but these improvements are quite minor for the present example. Pre-smoothing of the interface does not yield any improvement in the case of the monitor function~\eqref{eq:MackenzieMonitorFunction}; neither does pre-smoothing of the interface combined with smoothing on the subdomains, at least not in the measure $Q(K)$.

The above results confirm for the present example that the probabilistic approach to domain decomposition can yield good quality meshes at relatively cheap computational costs.   Different ways of smoothing have the potential to improve the mesh quality of domain decomposition grids further but ultimately may not be necessary at all.


\section{A Performance study}\label{sec:PerformanceStudy}

In this section we explore and comment on the computational performance of the proposed stochastic domain decomposition approach to grid generation. We consider the arc-length monitor function
\begin{subequations}\label{eq:FiveRingMonitorFunction}
\begin{equation}
 \rho=\sqrt{1+\alpha(u_x^2+u_y^2)},
\end{equation}
with velocity field of the form
\begin{align}
\begin{split}
 u&=\tanh\left[R\left(x^2+y^2-\frac18\right)\right]+
   \tanh\left[R\left(\left(x-\frac12\right)^2+\left(y-\frac12\right)^2-\frac18\right)\right]\\
   &+\tanh\left[R\left(\left(x-\frac12\right)^2+\left(y+\frac12\right)^2-\frac18\right)\right]+
     \tanh\left[R\left(\left(x+\frac12\right)^2+\left(y-\frac12\right)^2-\frac18\right)\right]\\
   &+\tanh\left[R\left(\left(x+\frac12\right)^2+\left(y+\frac12\right)^2-\frac18\right)\right],
\end{split}
\end{align}
\end{subequations}
on the domain $\Omega_{\rm p}=[-1,1]\times[-1,1]$, where $\alpha=0.2$ and $R=30$. Note that the same velocity field has been used in a similar monitor function in~\cite{huan05a}.

\begin{figure}[!ht]
\centering
\begin{subfigure}[b]{0.45\textwidth}
  \centering
  \includegraphics[width=\linewidth]{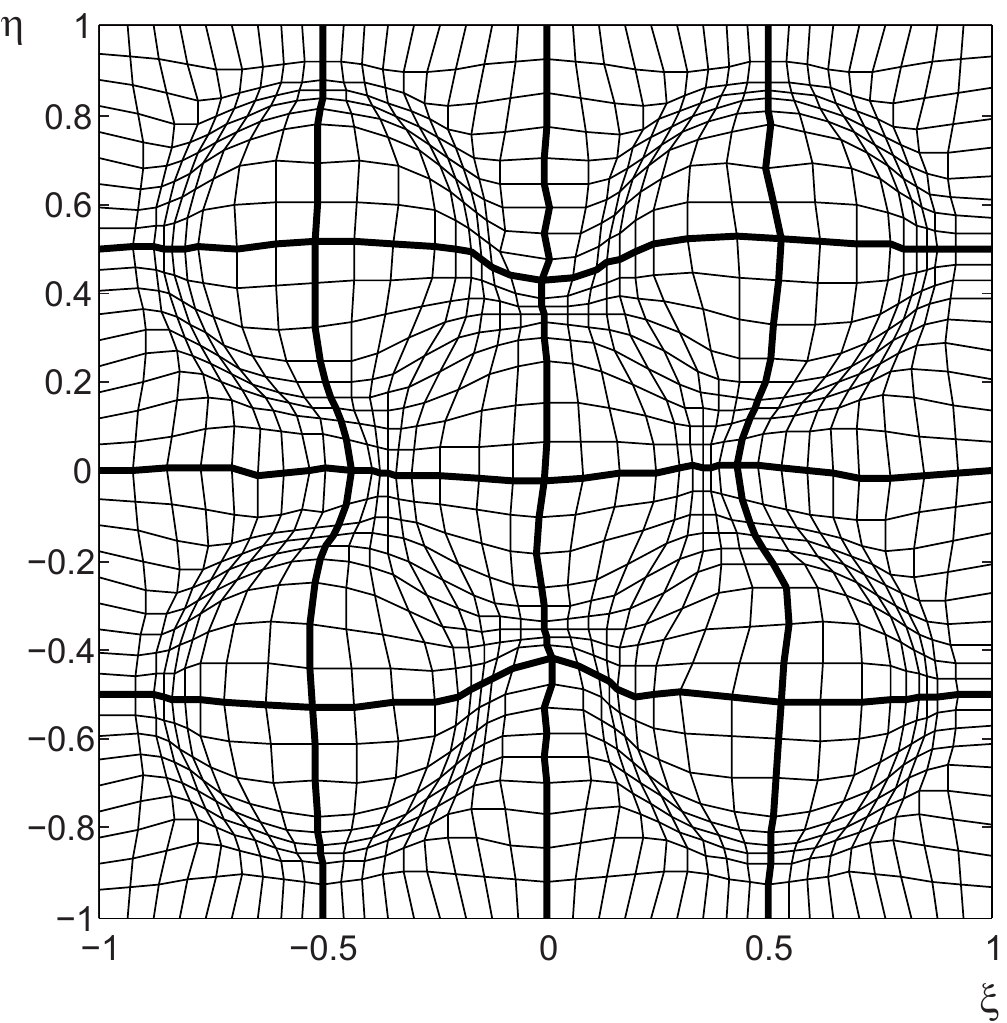}
  \label{fig:FiveRing410DD}
\end{subfigure}\qquad
\begin{subfigure}[b]{0.45\textwidth}
  \centering
  \includegraphics[width=\linewidth]{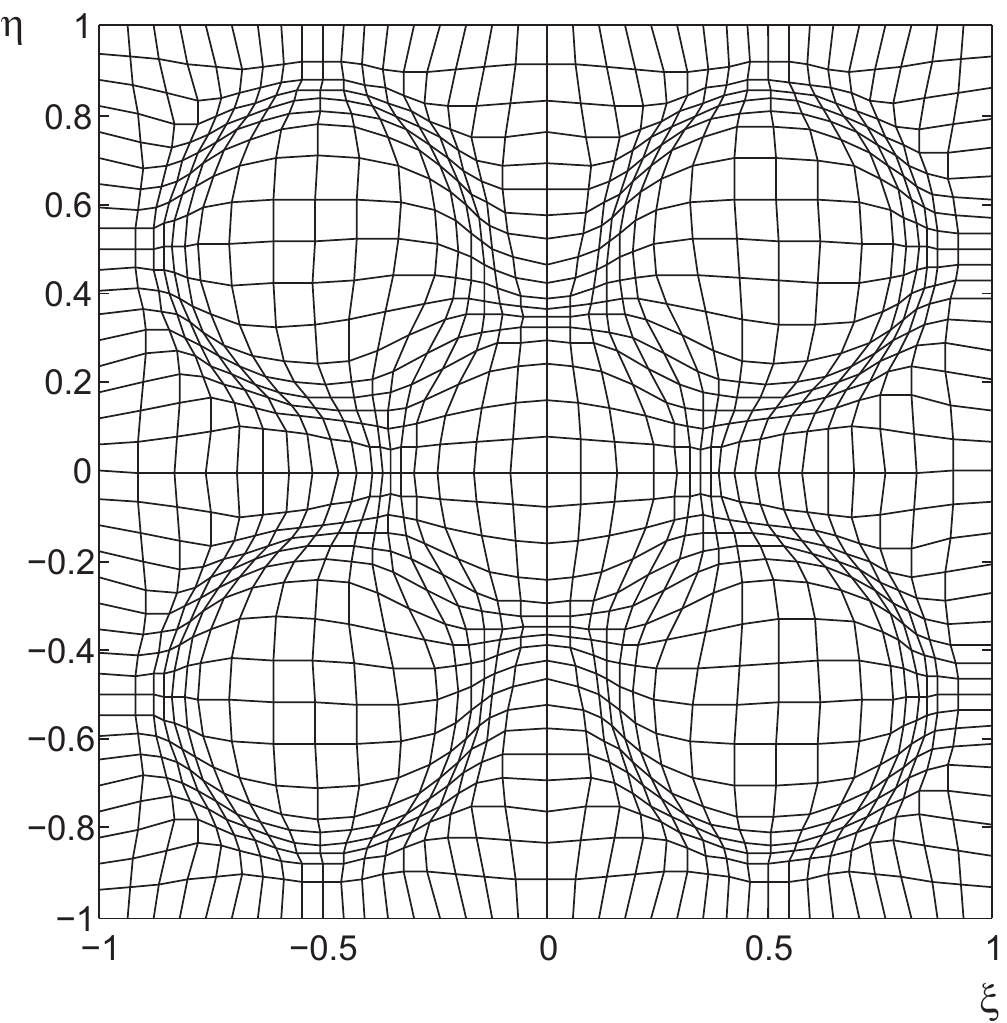}
  \label{fig:FiveRing410SD}
\end{subfigure}
\caption{\textbf{Left:} Domain decomposition solution of the mesh generator~\eqref{eq:WinslowTypeGridGenerator} using the monitor function~\eqref{eq:FiveRingMonitorFunction} to generate a grid on 16 subdomains with a total of $37^2$
mesh points.  \textbf{Right:} The associated global single-domain solution.}
\label{fig:FiveRing410}
\end{figure}

We solve this problem on 16 subdomains using a total of $37^2, 77^2, 117^2, 157^2$ and $197^2$ mesh points. The parameters for the stochastic solver were $N=5000$ and $\lambda=1000$. The result of the domain decomposition solution with $37^2$ mesh points and the associated single domain global solution are shown in Fig.~\ref{fig:FiveRing410}. The meshes for the other numbers of mesh points look qualitatively the same but become harder to depict due to the decrease in line spacings.

We compare run-times of the parallel algorithm on a 16 core parallel machine with the non-parallel single domain global solution. Results of this study are reported in Table~\ref{tab:ParallelPerformance} along with the associated geometric mesh quality measures.

\begin{table}[!ht]
\centering
\caption{Performance of the stochastic domain decomposition method for
grid generation using 16 subdomains with an increasing number of total mesh points.
Timings are given in seconds to two decimal places.}

\begin{tabular}{|r|r|r|r|r|r|r|}
\hline
Mesh Points  & $R_{\rm max}$ & $R_{\rm mean}$ & $T_{\rm stoc}$ & $T_{\rm sub}$ &
$T_{\rm total}$ & $T_{1}$ \\
\hline
$37 \times 37$ &   1  &  1     & 5.47 & 0.00 & 5.47 & 0.21\\
$77 \times 77$ & 0.97 &  1     & 7.56 & 0.02 & 7.58 & 1.67\\
$117\times 117$ & 0.98 &  1     & 8.81 & 0.08 & 8.89 & 7.04\\
$157\times 157$ & 0.94 &  1     & 8.88 & 0.21 & 9.09 & 20.62\\
$197\times 197$ & 0.95 &  0.99  & 9.06 & 0.47 & 9.53 & 47.80\\
\hline
\end{tabular}
\label{tab:ParallelPerformance}
\end{table}


In Table~\ref{tab:ParallelPerformance}, $T_{\rm stoc}$ is the time (in seconds) spent on computing the interface solutions using the stochastic representation of the solution, $T_{\rm sub}$ is the time required to deterministically obtain the mesh on all subdomains and $T_{\rm total} = T_{\rm stoc}+T_{\rm sub}$ is the total time of the SDD algorithm. The computation time needed for the global, non-parallel generation of the mesh is~$T_{1}$.  The same linear solver is used the global mesh and for each subdomain solve.

It can be seen from Table~\ref{tab:ParallelPerformance} that the stochastic domain decomposition method outperforms the deterministic, global mesh generator as soon as 16 subdomains are used with a total
of $157\times 157$ mesh points. Interestingly, the time spent on the probabilistic interface solution is almost constant when more than  $77\times 77$ mesh points are used. The reason for this is that the number of points on the interface where the stochastic solution is computed does not change when increasing the number of total points due to the optimal placing strategy employed. Also, the mesh quality is approximately the same for all cases reported in Table~\ref{tab:ParallelPerformance} and shows that the domain decomposition method produces meshes with almost identical quality as the single domain solution. The reason for the slight deterioration in mesh quality with an increased number of points in each subdomain is
that interpolation is used at more points along the interfaces. That is, the interpolation error becomes more relevant. No smoothing strategies have been applied to this problem.

The exact cross--over point at which the SDD algorithm becomes more efficient
depends on the number of Monte Carlo simulations used at each interface point, the number of interface points, the choice of the deterministic subdomain solver, the accuracy required and the number of compute cores available.
Communication costs have not been reported explicitly in this table.  Indeed such communication times could affect the
cross--over point for small numbers of mesh points but would not affect the conclusion for larger meshes.

As seen in the previous sections there are two possible strategies presented to improve mesh quality: additional
Monte Carlo simulations so that the interface values become more accurate or smoothing the computed solutions.  A natural question is what effect does additional Monte Carlo simulations and smoothing have on the performance of the algorithm?   We pursue this in Table \ref{tab:PerformanceWithoutSmoothing} and \ref{tab:PerformanceWithSmoothing}.    We time the various pieces of the algorithm with increasing numbers of Monte Carlo simulations.  In Table \ref{tab:PerformanceWithoutSmoothing} no smoothing is applied, while Table \ref{tab:PerformanceWithSmoothing} has been computed by
smoothing the interface values and the subdomain solutions as described in Section
\ref{sec:DDStochastic}.

As mentioned previously, increasing the number of Monte Carlo simulations increases mesh quality and of course increases the cost of the stochastic portion of the algorithm $T_{\rm stoc}$ and
the total computing time $T_{\rm total}$.
The timing for the subdomain solves is not affected by increasing $N$.
The time required for smoothing $T_{\rm smooth}$ does not depend on $N$.  It depends
only on the number of interface points and number of mesh points in each subdomain -- which
are both kept constant during the experiment.
For this example, the benefit of smoothing is clear, greatly improving mesh quality without
adding significantly to the cost.  Indeed using $N=2000$ Monte Carlo simulations and smoothing
is able to obtain a better mesh quality, in significantly less time,  than using $N=5000$ Monte Carlo simulation without smoothing.

\begin{table}[!ht]
\centering
\caption{Performance of the stochastic domain decomposition method for grid generation as a function of the number of Monte Carlo simulations without smoothing. A total of
$157\times 157$ points are used and the solution is computed with 16 subdomains.}
\begin{tabular}{|r|r|r|r|r|r|r|}
\hline
$N$ & $R_{\rm max}$ & $R_{\rm mean}$ & $T_{\rm stoc}$ & $T_{\rm sub}$ & $T_{\rm total}$ \\
\hline
500  & 0.04 &  0.96  & 0.91 & 0.21 & 1.12 \\
1000 & 0.24 &  0.98  & 1.88 & 0.21 & 2.09 \\
2000 & 0.71 &  1     & 3.68 & 0.21 & 3.8  \\
5000 & 0.93 &  1     & 9.07 & 0.21 & 9.28 \\
\hline
\end{tabular}
\label{tab:PerformanceWithoutSmoothing}
\end{table}

\begin{table}[!ht]
\centering
\caption{Performance of the stochastic domain decomposition method for grid generation as a function of the number of Monte Carlo simulations using
smoothing.  A total of
$157\times 157$ points are used and the solution is computed with 16 subdomains. }
\begin{tabular}{|r|r|r|r|r|r|r|}
\hline
$N$ & $R_{\rm max}$ & $R_{\rm mean}$ & $T_{\rm stoc}$ & $T_{\rm smooth}$ & $T_{\rm sub}$ & $T_{\rm total}$ \\
\hline
500  & 0.46 &  0.98  & 0.91 & 0.01 & 0.21 & 1.13 \\
1000 & 0.82 &  0.99  & 1.88 & 0.01 & 0.21 & 2.10 \\
2000 & 0.96 &  1     & 3.68 & 0.01 & 0.21 & 3.9  \\
5000 & 0.98 &  1     & 9.07 & 0.01 & 0.21 & 9.29 \\
\hline
\end{tabular}
\label{tab:PerformanceWithSmoothing}
\end{table}

In Table~\ref{tab:ParallelPerformanceConstantNumberOfPoints} we fix the problem size
and generate a $197\times 197$ mesh using a $2\times 2$, $3\times 3$ and $4\times 4$
configuration of subdomains.  As expected the time spent on each subdomain solve $T_{\rm sub}$
decreases as the number of subdomains increases.  The time spent on the stochastic portion
of the algorithm does not decrease dramatically with the number of subdomains as our
interface point placement strategy uses a constant number of mesh points on each interface
based on the local extrema and inflection points of the mesh density function.
This placement strategy, beneficial for larger subdomains, may be relaxed for larger numbers
of subdomains.

\begin{table}[!ht]
\centering
\caption{Scaling study for the generation of a $197\times 197$ mesh
 with a varying number of  subdomains $S$.
Timings are given in seconds to two decimal places.}
\begin{tabular}{|r|r|r|r|r|r|r|}
\hline
$S$ &  $R_{\rm max}$ &
$R_{\rm mean}$ & $T_{\rm stoc}$ & $T_{\rm sub}$ &
$T_{\rm total}$ & $T_{1}$ \\
\hline
$4$ &  0.96 & 1 & 12.54 & 14.89 & 27.43 & 64.00 \\
$9$ &  0.89 & 1 & 10.81 &  2.95 & 13.76 & 64.00 \\
$16$ & 0.96 & 1 & 10.59 &  0.46 & 11.05 & 64.00 \\
\hline
\end{tabular}
\label{tab:ParallelPerformanceConstantNumberOfPoints}
\end{table}

\section{Conclusions and outlook}\label{sec:ConclusionsDD}

In this paper we proposed a new method for computing PDE based adaptive meshes using a stochastic
domain decomposition approach.   The main idea derives from the study proposed in~\cite{aceb05a} to use stochastic domain decomposition to solve linear boundary value problems. As various types of adaptation strategies fit into this class of problems, this probabilistic take on domain decomposition is readily applicable to the mesh generation problem. The main motivation behind this work is to find an algorithm for the construction of meshes that is fully parallelizable as it allows one to determine the global mesh without information exchange across the interface of neighboring subdomains.
In addition, the relatively low accuracy requirement for mesh generation suggests that quality meshes can
be found with only a moderate number of Monte Carlo simulations along the interface or by
smoothing the obtained meshes;  hence
this mesh generation problem is well suited to this stochastic domain decomposition approach.

We have restricted ourselves to a simple linear mesh generator of Winslow type. Despite being amongst the simplest PDE based mesh generators, we found the Winslow potential system suitable to demonstrate several typical issues to be addressed when using a combined stochastic solver. These issues include ways to smooth the mesh in order to avoid an excessive number of Monte Carlo simulations as well as the optimal placement problem to determine the locations at the subdomain interfaces where the solution is needed probabilistically.  The former issue can be conveniently tackled using anisotropic diffusion to smooth out kinks that are inevitable due to the local nature of the Monte Carlo simulations, while properly preserving the regions of grid concentration and de-concentration. The latter problem of computing only few points along the subdomain interfaces and hence further reducing the number of Monte Carlo simulations needed is more difficult and ultimately boils down to a trade--off between the required grid smoothness and available computational resources. We found that placing the points near the maxima and minima of the first and second derivatives of the monitor function along the interfaces yields a sufficiently smooth mesh that can be determined at a fraction of the computational cost incurred if the solution at
all interface points are obtained probabilistically. More sophisticated criteria for this optimal placement problem will be investigated in the future.

A simple performance study shows the viability of the approach both as the problem size becomes larger and as the number of subdomains increases for a fixed problem size.
We also show that the cost of the stochastic portion of the algorithm can be kept small with the
incorporation of smoothing with little negative impact on the efficiency of the algorithm.

Theoretically, this approach is well--suited and theoretically extensible to three spatial dimensions as both
the Winslow linear mesh generator and the stochastic representation of the point--wise solution generalizes in a straightforward manner.  Moreover, the relative efficiency of the
stochastic approach increases as the spatial dimension of the problem increases.

A more challenging problem is the generalization of the proposed method to nonlinear mesh generators. It is well known that linear mesh generators of Winslow type work reasonably well for problems for which isotropic mesh generation suffices. In turn, if more control over the grid adaptation is required, one usually has to resort to nonlinear mesh generators, e.g.\ those that are based on equidistribution and alignment, see e.g.~\cite{huan01a,huan10a}. The problem is that for nonlinear partial differential equations, fewer stochastic solution representations are known. As the existence of a stochastic representation is at the heart of the proposed algorithm, it is crucial to either find a nonlinear mesh generator for which such a stochastic solution exists or to modify the problem so that it fits again into the realm of linear partial differential equations. The latter possibility may be realized through an appropriate linearization of the nonlinear mesh equations.
This is the subject of current investigations.

A further crucial step to refine the proposed algorithm
will be the development of a \emph{stopping criterion} for the Monte Carlo simulations. The task of generating meshes by solving a system of PDEs is different from the problem of obtaining numerical solutions of physical PDEs.  For the former problem the mesh PDE  is usually only a means to realize underlying computational grids with specific properties (e.g.\ isotropic vs.\ anisotropic adaptation). This means that it is usually not necessary to solve the mesh generation system with high accuracy, i.e.\ one can stop the solution procedure once a suitably smooth mesh is obtained. The derivation of proper stopping criteria for the mesh generation process will be another important source for improving the computational cost of the proposed algorithm. The main idea is to continuously evaluate the quality of the current candidate mesh and to stop or continue the solution procedure depending on whether a specified mesh quality threshold is reached. This strategy is fully compatible with the probabilistic computation of the interface solution due to the additive properties of the expected value.

As mentioned, the standard way to speed up Monte Carlo simulations is to use quasi-random numbers instead of pseudo-random numbers. Using quasi-random numbers was a crucial factor in~\cite{aceb05a} to obtain more accurate numerical solutions using the stochastic domain decomposition method. In this study we decided to use pseudo-random numbers only, as quasi-random numbers introduce a
sequential step in the algorithm.  Introducing this potential bottleneck in the algorithm does not appear justified for the mesh generation problem since a good quality mesh can be found without running the stochastic simulation to convergence. We reserve the more careful study of this issue for future work.

In addition, it is necessary to do a
careful comparison between this stochastic DD approach and the more
standard application of DD for mesh generation, cf.\ \cite{Haynes:2012, hayn12a}.   The candidates for comparison include an
optimized Schwarz approach, using DD as a preconditioner while solving
(\ref{eq:linearmesh})
in a Newton--Krylov--Schwarz method, and an ASPIN framework for nonlinear mesh generation.
See \cite{cai2002,cai2002b,cai2004,cai2005} for more details.

\section*{Acknowledgements}

The authors thank Weizhang Huang and the two anonymous referees for helpful remarks regarding the manuscript. This research was supported by the Austrian Science Fund (FWF), project J3182--N13 and the Natural Sciences and Engineering Research Council of Canada (NSERC). AB is a recipient of an APART Fellowship of the Austrian Academy of Sciences.

{\footnotesize\itemsep=0ex
\def\cprime{$'$}

}

\end{document}